\documentclass[twoside]{article}

\usepackage[accepted]{aistats2021}

\usepackage[utf8]{inputenc} 
\usepackage[T1]{fontenc}    
\usepackage{hyperref}       
\usepackage{url}            
\usepackage{booktabs}       
\usepackage{amsfonts}       
\usepackage{nicefrac}       
\usepackage{microtype}      
\usepackage{bm}
\usepackage[round]{natbib} 
\usepackage{xfrac}
\usepackage{comment}
\usepackage{amsthm}
\usepackage{thm-restate}
\newtheorem{thm}{Theorem}
\newtheorem{lem}[thm]{Lemma}
\newtheorem{cor}[thm]{Corollary}

\usepackage{enumitem}
\usepackage{tcolorbox}
\tcbset{boxsep=0mm,boxrule=0pt,colframe=white,arc=0mm,left=0.5mm,right=0.5mm}

\usepackage{amsmath,thmtools}

\DeclareMathOperator*{\argmin}{arg\,min}
\renewcommand{\rm}{{\mathcal R}}

\newcommand{\R}{{\mathbb R}}

\newcommand{\sie}{\text{\tiny SIE}}
\newcommand{\ee}{\text{\tiny EE}}

\theoremstyle{remark}
\newtheorem{remark}{Remark}

\begin{document}

\runningauthor{Zhang, Orvieto, Daneshmand, Hofmann, Smith}
\runningtitle{Revisiting the Role of Symplectic Integration on Acceleration and Stability in Convex Optimization}

\twocolumn[
\aistatstitle{Revisiting the Role of Euler Numerical Integration\\ on Acceleration and Stability in Convex Optimization}

\aistatsauthor{Peiyuan Zhang* \And Antonio Orvieto \And Hadi Daneshmand \And Thomas Hofmann \And Roy Smith}

\aistatsaddress{ETH Z\"urich  \And ETH Z\"urich \And INRIA Paris, ETH Z\"urich \And ETH Z\"urich \And ETH Z\"urich } 
]

\begin{abstract}
Viewing optimization methods as numerical integrators for ordinary differential equations (ODEs) provides a thought-provoking modern framework for studying accelerated first-order optimizers. In this literature, acceleration is often supposed to be linked to the quality of the integrator (accuracy, energy preservation, symplecticity). In this work, we propose a novel ordinary differential equation that questions this connection: both the explicit and the semi-implicit~(a.k.a symplectic) Euler discretizations on this ODE lead to an accelerated algorithm for convex programming. Although semi-implicit methods are well-known in numerical analysis to enjoy many desirable features for the integration of physical systems, our findings show that these properties do not necessarily relate to acceleration.
\end{abstract}

{\let\thefootnote\relax\footnotetext{*Correspondence to talantyeri@gmail.com.}}

\section{Introduction}
Momentum methods are the state-of-the-art choice of practitioners for the optimization of machine learning models. The simplest of such algorithms is the Heavy-ball~(HB), first proposed and analyzed in the context of convex optimization by \citet{polyak1964some}:
\begin{equation}
    \label{eq: hb} \tag{HB}  x_{k+1} = x_k + \beta(x_k-x_{k-1}) - s \nabla f(x_k)
\end{equation}
where $f:\R^d\to\R$ is the $L$-smooth\footnote{For all $x,y\in\R^d$, $\|\nabla f(x)-\nabla f(y)\|\le L\|x-y\|$, where $\|\cdot\|$ is the standard Euclidean norm.} function we want to minimize, $s>0$ is the step-size and $\beta\in[0,1)$ the momentum parameter. Using a novel and beautiful argument on fixed point iterations, \citet{polyak1964some} proved that, if $f$ is twice continuously differentiable and $\mu$-strongly-convex\footnote{ $\forall x\in\R^d, $ $\nabla^2 f(x)-\mu I$ is positive semidefinite.}, the sequence $(x_k)_{k\ge0}$ produced by \ref{eq: hb} \textit{locally}~(i.e. if initialized close to the solution) converges to the minimizer $x^*=\argmin_{x\in\R^d} f(x)$ at an \textit{accelerated} rate. The keyword ``accelerated'' has a precise meaning: an algorithm for $\mu$-strongly-convex and $L$-smooth problems is accelerated if and only if the convergence rate of $f(x_k)$ to $f^* := \min_{x\in\R^d} f(x)$ is $O((1-\sqrt{\mu/L})^k)$. For instance, Gradient Descent~(i.e. $\beta=0$) in this setting converges linearly but with constant $1-\mu/L$ and is therefore not accelerated\footnote{If $L/\mu$ is large, $1-\sqrt{\mu/L}\ll 1-\mu/L$.}.
\vspace{-2mm}
\paragraph{Nesterov's acceleration.} Supported by the lower bounds established by \citet{nemirovsky1983problem}, many researchers in the early 80s tried to develop an algorithm with a \textit{global} accelerated convergence rate. The problem was solved by \citet{nesterov1983method}, who proposed the following modification\footnote{In the original paper~\citet{nesterov1983method}, the algorithm is presented in a more general way. Our formulation is similar to~\citet{shi2018understanding}.} of \ref{eq: hb}:
\begin{align}
	\label{eq: nag} \tag{NAG} 
	\begin{split} x_{k+1} = x_k &+ \beta(x_k-x_{k-1}) - s \nabla f(x_k) \\
	&- \beta s(\nabla f(x_k) - \nabla f(x_{k-1})).
	\end{split}
\end{align}
The intuition behind this algorithm puzzled researchers for decades, and many articles are devoted to understanding the underlying mechanism~\citep{ allen2014linear,defazio2019curved,ahn2020proximal} and the role of the small yet crucial modification\footnote{This is usually referred to as \textit{gradient extrapolation}.} compared to \ref{eq: hb}~\citep{flammarion2015averaging,lessard2016analysis,hu2017dissipativity}. Notwithstanding the theoretical value of these contributions, they are arguably only of a descriptive nature and leave open more fundamental questions on the reason behind acceleration.
\vspace{-2mm}
\paragraph{Continuous-time models for acceleration.} A new line of research bloomed from a seminal paper by \citet{su2014differential}. This work gained a lot of attraction, as it introduces\footnote{We point out that, actually, the differential equations proposed in~\citet{su2014differential} was already written down and partly analyzed in the original 1963 paper by \citet{polyak1964some}. Even more surprisingly, a first study of damped second order differential equations for optimization can be found already in a 1958 paper of the soviet mathematician \citet{gavurin1958nonlinear}. } a powerful way to look at acceleration through the lens of second order ordinary differential equations~(ODEs). In the $\mu$-strongly-convex case, this equation is 
\begin{align}
	\label{eq: nag_ode} \tag{NAG-ODE} 
\ddot{X} +2\sqrt{\mu}\dot X +\nabla f(X)=0
\end{align}
and retains the essence of acceleration: namely, convergence with a rate $O(e^{-\sqrt{\mu} t})$. Analogously to the discrete-time case we just discussed, one can prove that the continuous-time model of gradient descent, i.e. the gradient flow $\dot X = -\nabla f(X)$, converges instead at the non-accelerated rate $O(e^{-\mu t})$. Other interesting properties of damped gradient systems such as \ref{eq: nag_ode} can be found in the (stochastic) optimization literature~\citep{krichene2015accelerated,xu2018continuous,cabot2009long,orvieto2020role,orvieto2019shadowing,diakonikolas2019generalized,alecsa2019extension,alimisis2020continuous}, and in the applied mathematics literature~\citep{sanz2020connections,attouch2000heavy,attouch2000heavy2,alvarez2000minimizing,begout2015damped}.
\vspace{-2mm}
\paragraph{High-resolution ODEs.} As first noted by~\citet{wilson2016lyapunov}, while \ref{eq: nag_ode} is formally the continuous-time limit~(for some specific choice of $\beta$) of \ref{eq: nag}, it is also the continuous-time limit of \ref{eq: hb}. In other words, \ref{eq: nag_ode} does not capture the vanishing gradient correction~(a.k.a gradient extrapolation) term $\beta s(\nabla f(x_k) - \nabla f(x_{k-1}))$, which is regarded to be a fundamental piece of the acceleration machinery in discrete-time. To solve this issue~(i.e. to get a more accurate model of Nesterov's acceleration), \citet{shi2018understanding} introduced a high-resolution model of \ref{eq: nag}:
\begin{align}
	\label{eq: nag-hr}\tag{NAG-ODE-HR}
	\begin{split}
\ddot{X} & + (2\sqrt{\mu}  + \sqrt{s}\nabla^2 f(X))\dot{X} \\
&+ (1 + 2\sqrt{\mu s})\nabla f(X) = 0.
\end{split}
\end{align}
Remarkably, here (1) the step-size $s$ is included directly in the model, and (2) the vanishing (as $s\to0$) term $\sqrt{s}\nabla^2 f(X)\dot X$ is used to capture the gradient correction $\beta s(\nabla f(x_k) - \nabla f(x_{k-1}))$. The term $\sqrt{s}\nabla^2 f(X)\dot X$ is referred to as Hessian damping, and can be seen as a curvature-dependent viscosity correction. As a validation for their new ODE, \citet{shi2018understanding} showed that \ref{eq: nag-hr} enjoys the same accelerated rate of \ref{eq: nag_ode} --- but it is empirically more faithful to \ref{eq: nag} compared to \ref{eq: nag_ode}, for finite values of $s$.
\vspace{-2mm}
\paragraph{Connection to numerical integration.} In a second article, \citet{shi2019acceleration} showed that~\ref{eq: nag} can be approximately recovered\footnote{See approximations in Sec.~2.2 of~\citet{shi2019acceleration}.} through a semi-implicit~(a.k.a. symplectic) Euler discretization of \ref{eq: nag-hr}. The authors also claim that if the same system is integrated with the explicit Euler method, the resulting optimizer \textit{might}\footnote{We point out here a potential problem in the main claim of~\cite{shi2019acceleration}: the authors show that the explicit integrator of \ref{eq: nag-hr} is stable only for small step-sizes by finding necessary conditions for the steady decrease of a particular energy function. While this fact surely hints at potential instabilities of the associated algorithm, it does not \textit{per se} provide a sufficient condition for slow convergence. In other words, one could in principle find a different result by choosing a different Lyapunov function.} not be accelerated because it is only found stable for small values of $s$. Semi-implicit methods are well-known to perform remarkably well for integrating second-order ODEs in physics~\citep{hairer2006geometric} and chemistry~\citep{lubich2008quantum}; namely, one can use big step-sizes while preserving the geometry of the original flow. \textit{Shi et al. suggested that the essence of acceleration can be explained by the same phenomenon}, which is mathematically well understood in the Hamiltonian~(i.e. energy-conserving) setting thanks to the theory of backward error analysis~\citep{hairer1994backward,benettin1994hamiltonian}.

\begin{figure}[ht]
    \centering
 	\includegraphics*[width=0.75\linewidth]{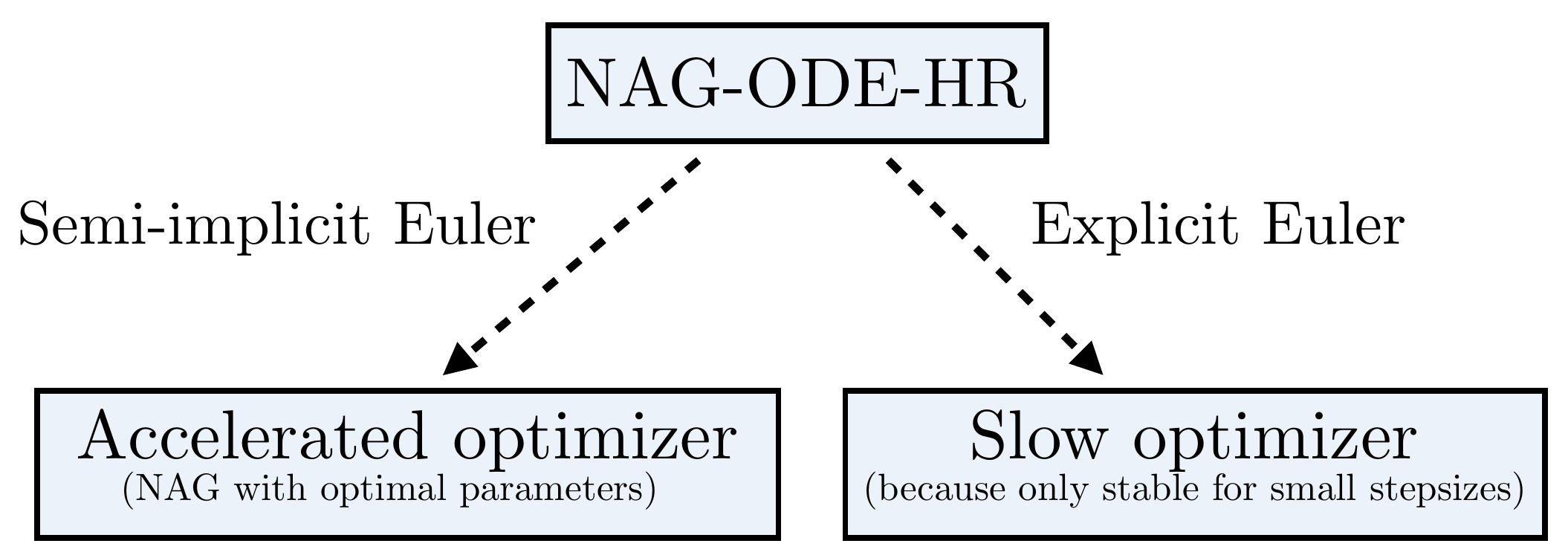} 
 	\caption{\small Sketch of the storyline of \citet{shi2019acceleration}: while semi-implicit discretization of \ref{eq: nag-hr} yields an accelerated method, explicit discretization results in a method not known in the literature, which is claimed to be stable only for very small step-sizes ($s\le O(\mu/L^2)$, compared to $s\le O(1/L)$ of the semi-implicit method). This is used by the authors to advocate that the structure provided by semi-implicit integration is somethow critical for the construction of accelerated methods. This storyline~(and the associated conclusion) is much different from ours, sketched in Fig.~\ref{fig:our_picture}. Superiority of semi-implicit methods is also claimed/hinted in several works~\citep{shi2018understanding,bravetti2019optimization,betancourt2018symplectic,francca2020dissipative,muehlebach2019dynamical}.}
 	\label{fig:shi_picture}
\end{figure}

On a parallel line, \citet{muehlebach2019dynamical} derived a different continuous-time model that contains terms of the form $\nabla f(X+\sqrt{s}\dot X)$ instead of $\sqrt{s}\nabla^2f(X) \dot X$. This ODE can also relate to Nesterov's method through semi-implicit integration. Moreover, inspired by the variational perspective presented in~\citet{wibisono2016variational}, many research papers~\citep{betancourt2018symplectic,muehlebach2020optimization,francca2020dissipative,francca2020conformal,alecsa2020long,bravetti2019optimization} have been devoted to understanding the geometric properties of Nesterov's method, seen as either (1) a (Strang/Lie-Trotter) splitting scheme for structure-preserving integration of conformal Hamiltonian systems~\citep{mclachlan2001conformal,mclachlan2002splitting} or (2) the composition of a map derived from a contact Hamiltonian~\citep{de2019contact,bravetti2017contact} and a gradient descent step. Finally, the application of Runge-Kutta schemes was explored \citep{zhang2018direct,zhang2019acceleration,sanz2020contractivity}; in particular,~\citet{zhang2018direct} first showed that fast rates can be also achieved via high-order explicit methods.

To sum it up, to the best of our knowledge, \textit{all} recent convex optimization literature advocates that, in order to achieve acceleration from an ODE model, one needs to use either structure-preserving integrators~\citep{bravetti2019optimization,shi2018understanding,francca2020dissipative}, high-order explicit methods~\citep{zhang2018direct,zhang2019acceleration}, or implicit methods~\citep{shi2019acceleration,wilson2016lyapunov,diakonikolas2017accelerated}.
\vspace{-2mm}

\paragraph{Our contribution.} We show that, contrary to what is often claimed~(or hinted at) in recent literature~(see paragraph above) acceleration can also be achieved by means of simple low-order explicit numerical integrators --- such as the explicit Euler method. While explicit Euler is well-known to be provably suboptimal for accurate integration of Hamiltonian systems \citep{hairer2006geometric,hairer1994backward,benettin1994hamiltonian}, we show that this does not necessarily imply slow convergence of the resulting optimizer. In particular, our work suggests that the structure-preserving properties of semi-implicit~(symplectic) methods are not a necessary component of accelerated algorithms. 

We start by introducing a generalized momentum ODE (GM-ODE), dependent on three parameters, which recovers both \ref{eq: nag_ode} and \ref{eq: nag-hr} as special cases. In Sec.~\ref{sec:cont_model} we study the convergence of this ODE. Next, in Sec.~\ref{sec:discretization} we show that both the explicit and the semi-implicit Euler methods, when applied for numerical integration of GM-ODE, can achieve an accelerated rate. Finally, in Sec.~\ref{sec:truncation}, we go one step further and show that there exist damped gradient systems for which the semi-implicit Euler method is unstable, while the explicit Euler method (with the same step-size) is stable. Of course, for other ODE systems, we observe the opposite behavior. This showcases that the stability of the integrator depends on the underlying accelerated ODE.

\begin{figure}[ht]
    \centering
 	\includegraphics*[width=0.75\linewidth]{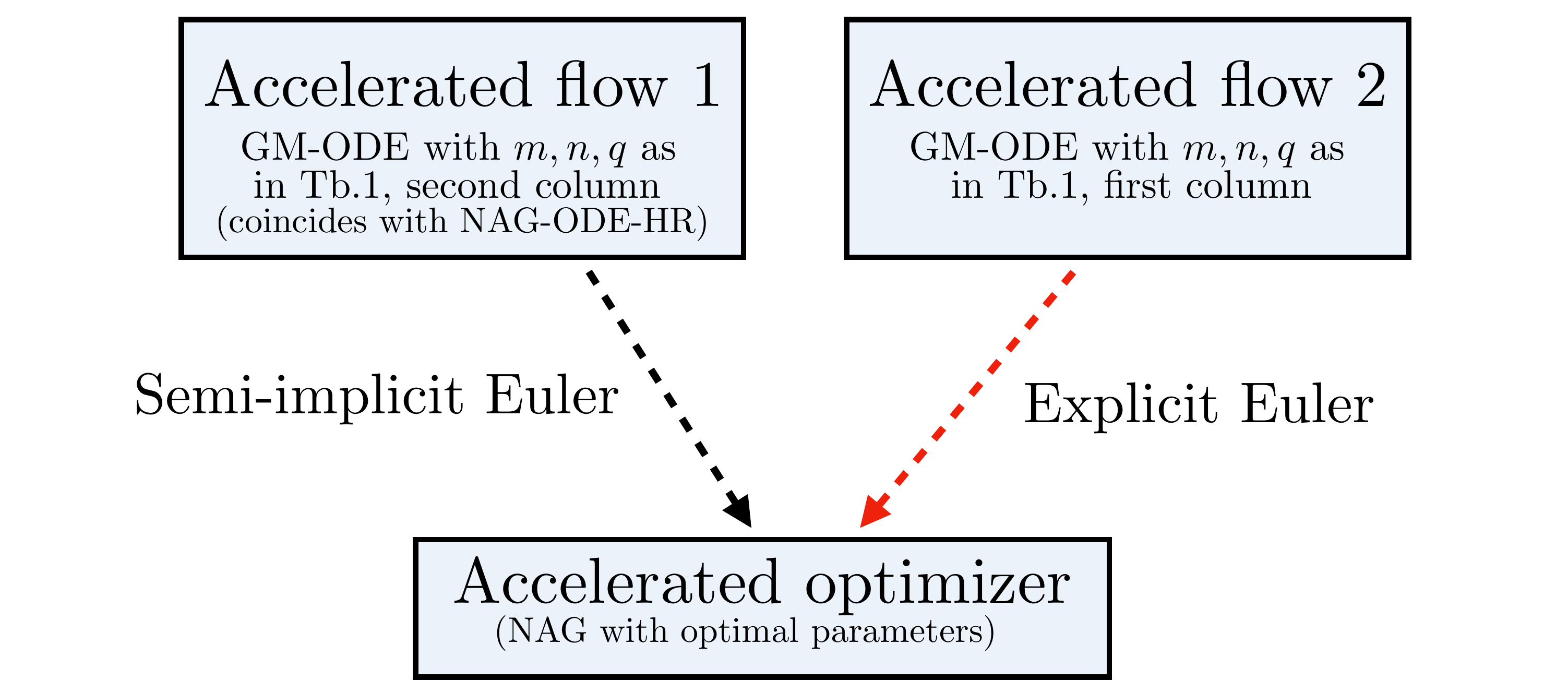} 
 	\caption{\small The main conceptual finding of our paper: both the semi-implicit and the explicit Euler integrators are able to recover Nesterov's method, when used to discretize different accelerated flows~(ODEs are given in Sec.~\ref{sec:discretization}).}
 	 	\label{fig:our_picture}
\end{figure}
\vspace{-2mm}
At its core, our work showcases some unintuitive aspects of the connection between the fields of numerical integration and optimization. Namely, while for accurate integration of physical systems symplectic integrators are provably superior to explicit methods~\citep{hairer2006geometric,benettin1994hamiltonian}, we show that the same ranking might not hold when seeking fast optimizers through ODE discretization. We think that clarifying this critical point makes the (already vast) literature on this topic richer, motivating future research on the connections between the fields of optimization and numerical analysis.

\section{Summary of the results}
\vspace{-2mm}
Our work is based on the study of a novel continuous-time model for momentum algorithms, namely the following ordinary differential equation that is indexed by non-negative parameters $m, n, q$:
\begin{align}
\begin{cases}
\dot{X} = - m \nabla f(X) - nV \\
\dot{V} = \nabla f(X) - qV,  \label{eq:our_model} \tag{GM-ODE}
\end{cases}
\end{align}
We show that the above ODE includes both \ref{eq: nag_ode} and \ref{eq: nag-hr} as special cases and recovers a large set of momentum methods through the application of two classical numerical integrators, i.e. semi-implicit and explicit Euler. In Lemma~\ref{lemma:equivalence} we show that these integrators are equivalent and can both lead to acceleration~(see also Fig.~\ref{fig:our_picture}). Equivalence is shown by reparameterization: any semi-implicit discretization of \ref{eq:our_model} with parameters $(m_{\sie},n_{\sie},q_{\sie})$ can be viewed as explicit Euler discretization of \ref{eq:our_model} with different parameters $(m_{\ee},n_{\ee},q_{\ee})$. These parameters can be computed in closed form starting from $(m_{\sie},n_{\sie},q_{\sie})$~(see Lemma~\ref{lemma:equivalence} for precise formulas). Such an equivalence suggests that the energy-preservation properties of semi-implicit integration might no be strictly necessary to achieve acceleration, as instead hinted in recent works~\citep{bravetti2017contact,mclachlan2001conformal,shi2019acceleration}.

To make our analysis complete, we establish an accelerated convergence rate (in Thm.~\ref{thm:convergence} and Cor.~\ref{thm:convergence_ee}) for a set of algorithms that can be interpreted as~(both the two) Euler discretizations of \ref{eq:our_model} with a proper parameters choice. 
As a side product of our novel analysis of semi-implicit and explicit methods, we derive a novel accelerated convergence rate for the quasi-hyperbolic momentum~(QHM) method introduced by \citet{ma2018quasi}. Indeed, along with \ref{eq: hb} and \ref{eq: nag}, the QHM method can also be seen as a numerical integrator on~\ref{eq:our_model}. QHM was shown to be very competitive in deep learning tasks~\citep{choi2019empirical} as well as in the strongly-convex setting~(see Appendix J in~\citep{ma2018quasi}). However, to the best of our knowledge, QHM has only been studied in the quadratic case~\citep{gitman2019understanding} (hence the novelty of our rate). We like to point out that this is not the main contribution of our paper, but is presented here nonetheless to showcases the flexibility of our novel ODE and of the numerical integration approach.

Finally, we go beyond convergence analysis and study the discretization errors in Sec.~\ref{sec:truncation}. Under some conditions on the choice of parameters, we show that the explicit Euler method enjoys the same integration error as the semi-implicit Euler method when integrating \ref{eq:our_model}~(see Lemma~\ref{lem:global_disc_bound}).   
 



\section{Continuous-time analysis}
\label{sec:cont_model}
Before discussing numerical integration, we provide here a continuous-time analysis of \ref{eq:our_model}, in line with most related works on acceleration and numerical integration~\citep{shi2018understanding,su2014differential}. The results in this section are not fundamental for the understanding of our claims on the discretization of \ref{eq:our_model}. Hence, for a quick read, this section can be safely skipped.

\ref{eq:our_model} can be seen as as a linear combination of the gradient flow $\dot{X} = -\nabla f(X)$ (obtained for $n=0$) and \ref{eq: nag_ode}~(obtained for $n=1$). Assuming the objective function $f$ is $L$-smooth, one can check that \ref{eq:our_model} admits a unique solution~(follows directly from Thm. 3.2 in~\citet{khalil2002nonlinear}). The model above is inspired by the quasi-hyperbolic momentum~(QHM) algorithm\footnote{QHM was introduced as weighted average of momentum and gradient descent methods. It is shown to recover both \ref{eq: hb} and \ref{eq: nag} as special cases.} developed in \citet{ma2018quasi}. We discuss the connection to QHM later in Sec.~\ref{sec:discretization}.
\vspace{-2mm}
\paragraph{Connections to existing ODE models.} \ref{eq:our_model} recovers existing continuous momentum models under different choices of parameters. To see this, let us take the second derivative of $X$: $\ddot{X} = - m \nabla^2 f(X) \dot{X} - n \dot{V}$.
\vspace{-1mm}
\begin{align}
\ddot{X} + \big(q + m\nabla^2 f(X)\big)\dot{X} + (n + qm)\nabla f(X) = 0.
\label{eq:QHM_to_jordan}
\end{align}
The choice\footnote{Proofs for discretized \ref{eq:our_model} will rely on condition $m>0$. This discussion will be elaborated in Sec.~\ref{sec:discretization}.} $m=0, n=1, q=2\sqrt{\mu}$ recovers \ref{eq: nag_ode} by \citet{polyak1964some}. Moreover, the choice $m=\sqrt{s}, n=1, q=2\sqrt{\mu}$ recovers~\ref{eq: nag-hr}, proposed by \citet{shi2018understanding,shi2019acceleration}. That is, \ref{eq:our_model} includes as special cases both the high-resolution and low-resolution models of Nesterov's method~(see discussion in the introduction). We note that, contrary to~\citet{shi2018understanding}, the Hessian of $f$ is not explicitly included in the model. Also, contrary to~\citet{muehlebach2019dynamical}, the gradient is evaluated only at the current position $X$. This feature arguably gives \ref{eq:our_model} higher interpretability than existing models -- a simple linear combination of gradient and momentum can also achieve high resolution\footnote{That is, a finer, compared to the original ODE in~\citet{su2014differential} approximation of Nesterov's method. For a detailed discussion on this terminology, we refer the reader to~\citet{shi2018understanding}.}.

\begin{center}
\begin{tabular}{ |c||c|c|c| } 
 \hline
 &$m$ & $n$ & $q$ \\ 
\hline\hline
 Gradient Flow & 1 & $0$ & any \\
  \hline
 \ref{eq: nag_ode}~\citep{su2014differential}& $0$ & $1$ & $2\sqrt{\mu}$ \\ 
 \hline
 \ref{eq: nag-hr}~\citep{shi2019acceleration}& $\sqrt{s}$ & $1$ & $2\sqrt{\mu}$ \\
    \hline
\end{tabular}
\end{center}

\vspace{-2mm}
\paragraph{Stability and convergence rate.} The equilibria of~\ref{eq:our_model} are easy to characterize: since $m,n$ and $q$ are non-negative, we have $\dot X=0$ and $\dot V=0$ if and only if both $\nabla f(X)=0$ and $V=0$. Under the assumption that $f$ is strongly-convex, only its unique minimizer $x^*$ is such that $\nabla f(x^*)=0$. Therefore the point $(x^*,0)\in\R^{2d}$ is the only equilibrium of~\ref{eq:our_model}. Next, we want to show that $(x^*,0)$ is asymptotically stable and characterize the convergence rate of our model. Borrowing some inspiration from~\citet{su2014differential,shi2019acceleration}, we propose the following energy function:
\begin{align*}
\begin{split}
	\mathcal{E}(X,V) &= (qm+n)\big(f(X)-f(x^*)\big) \\
	&+ \frac{1}{4}\|q(X-x^*)-nV\|^2 + \frac{n(qm+n)}{4}\|V\|^2.
\end{split}
\end{align*}
The next theorem states our result about Lyapunov stability, of which the proof is presented in the appendix. 
\begin{tcolorbox}
\begin{restatable}[Continuous-time stability]{thm}{continuous}\label{thm:continuous}
	Let $f$ be $\mu$-strongly-convex and $L$-smooth. If $n,m,q\ge0$ then, for any value of the strong-convexity modulus $\mu\ge0$, the point $(x^*,0)\in\R^{2d}$ is globally asymptotically stable for~\ref{eq:our_model}, as
	\begin{align}
	\mathcal{E}(X(t),V(t)) & \leq e^{-\gamma_1 t} \cdot \mathcal{E}(X(0),V(0)),
	\label{eq:rate_conti}
	\end{align}
	where $\gamma_1 := \min\left(\cfrac{\mu(n+qm)}{2q}, \cfrac{q}{2}\right)$.
\end{restatable}
\end{tcolorbox}

\begin{figure}
    \centering
    \includegraphics[width=\linewidth]{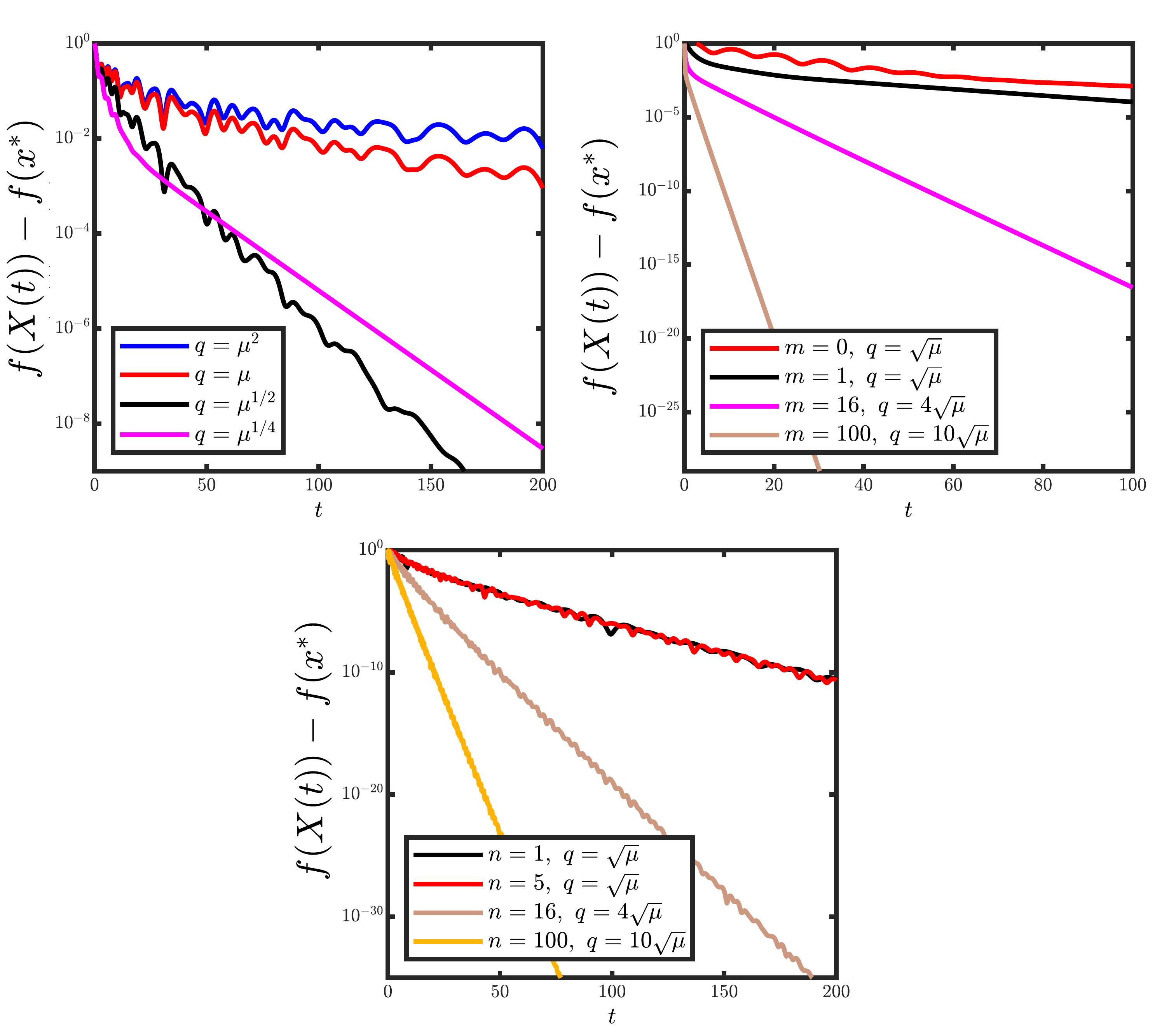}
\caption{\small{Role of parameters in \ref{eq:our_model}. The objective function is a 10-dimensional quadratic function $f$ with $\mu I\preceq\nabla^2 f\preceq L I$ where $\mu=0.01$ and $L=1$. The panels depict, from left to right, the influence of $q$, $m$ and $n$, as suggested in above discussion. In each figure we vary the parameter we are interested in (as in the legends) and keep the others fixed. For left panel we use $m=0.2$ and $n=1$. For the middle panel we use $n=0.1$. For the right panel we use $m=0.2$. Numerical integration of \ref{eq:our_model} performed using a fourth-order Runge-Kutta with step-size $10^{-4}$. }}
	\label{fig:continuous}
\end{figure}

Remarkably, the stability analysis in the proof can be used to guide the analysis of different momentum methods~(see Sec.~\ref{sec:discretization}) --- obtained by the application of standard Euler integrators of our model.
\begin{remark}
The rate in Thm.~\ref{thm:continuous} is not affected by the gradient Lipshitz constant $L$. This might look strange at first for a reader familiar with the optimization literature. However, we point to the fact that this is a feature of most continuous-time models~(see e.g. rates in~\citet{shi2018understanding}). The Lipschitz constant comes back into the rate after discretization, since one has to introduce a bound on the maximum integrator step-size, usually proportional to $1/L$~(see Eq.~\ref{eq:parameters_sie} and \ref{eq:parameters_ee}). 
\end{remark}
\vspace{-2mm}
\paragraph{How do $\boldsymbol{m,n,q}$ affect the ODE dynamics?} One can readily check that Thm.~\ref{thm:continuous} implies a linear rate in function value of the form $f(X(t))-f(x^*) \le O(-e^{\gamma_1 t})$. This result recovers exactly the rates in~\citet{shi2018understanding} as a special case. However, we note that our result is more general and leads to novel insights on the interplay between gradient amplification~(controlled by $n$), momentum~(controlled by $q$) and Hessian damping~(controlled by $m$). Indeed, given the expression for $\gamma_1$ in Eq.~\ref{eq:rate_conti}, we can make the following conclusions.
\setlist[2]{noitemsep} 
\setenumerate{noitemsep} 
\begin{itemize}[noitemsep,leftmargin=*]
\vspace{-3mm}
    \item For fixed $m,n\ge0$, the value of $q$ which maximizes $\gamma_1$ also solves $\mu(n+qm)/q = q$, which implies $q=(m+\sqrt{4\mu n + m^2})/2$. If we restrict $q$ to be a power of $\mu$, set $n=1$ and ignore the effect of $m$, then we get the popular choice~\citep{shi2018understanding,shi2019acceleration,muehlebach2019dynamical} $q = O(\sqrt{\mu})$~(see the first panel of Fig.~\ref{fig:continuous}). \citet{sanz2020connections} recently showed that this choice is optimal using the linear matrix inequalities framework~\citep{lessard2016analysis,fazlyab2018analysis}.
    \item For any $q\ge0$, if $n\ge0$ is chosen small enough such that $q^2-\mu n\ge 0$, then by picking $m = (q^2-\mu n)/q$ we have $\gamma_1 = q/2$. Hence, by increasing $q$ (and adapting $m$ accordingly) the convergence in continuous-time can be sped-up arbitrarily~(see the second panel of Fig.~\ref{fig:continuous}).
    \item If $n = q^2/\mu$, then $\gamma_1 = q/2$ for all $q\ge0$ and any $m\ge0$. Again, by increasing $q$ the convergence can be sped-up arbitrarily~(bottom panel of Fig.~\ref{fig:continuous}).
\end{itemize}
\begin{remark}
If $n$ or $m$ are increased, one can guarantee arbitrarily fast convergence to the minimizer. This result only holds true in continuous-time, as noted also in a similar setting by~\citet{wilson2016lyapunov}. Indeed, as we will see in Thm.~\ref{thm:convergence}, in the discrete word, \textit{to ensure stability}, $n$ and $m$ have to be bounded by a constant which is inversely proportional to the discretization accuracy.
\label{rk:accuracy_stability}
\end{remark}

\section{Discretization and acceleration}
\label{sec:discretization}
We now jump to the discrete world and show how \textit{both} explicit and semi-implicit numerical integration, applied to~\ref{eq:our_model}, can yield accelerated gradient iterations.
\vspace{-2mm}
\paragraph{Discretization schemes.} We consider two well-understood~\citep{hairer2006geometric} and practical first-order numerical integration schemes applied to~\ref{eq:our_model} with discretization step-size $\sqrt{s}$~(see discussion in~\citet{su2014differential,shi2019acceleration}): Explicit Euler~(EE) and Semi-Implicit\footnote{Actually, there exist many semi-implicit methods that go under the name of ``semi-implicit Euler''. We expect many of those integrators to work equally well for the sake of our discussion on equivalence. For a more detailed discussion, we refer the reader to Chapter 1 of~\citet{hairer2006geometric}.} Euler~(SIE).  
\allowdisplaybreaks
\begin{align*}
  \text{(EE)} : & \begin{cases}
x_{k+1} - x_k = - m\sqrt{s} \nabla f(x_k) -  n \sqrt{s} v_k \\
v_{k+1} - v_k = \sqrt{s} \nabla f(x_{k}) - q \sqrt{s} v_{k}.
\end{cases}\\ 
   \text{(SIE)} : &\begin{cases}
x_{k+1} - x_k = - m\sqrt{s} \nabla f(x_k) - n \sqrt{s} v_k \\
v_{k+1} - v_k = \sqrt{s} \nabla f(x_{k+1}) - q\sqrt{s} v_{k}.
\end{cases}
\end{align*}
Even though the second equation in SIE is written in an \textit{implicit} way, it can be trivially solved: indeed, one shall first find $x_{k+1}$ and then plug the solution into the second equation. Since the gradient computed at $x_{k+1}$ can be used at the next iteration, the two algorithms have the same complexity. Indeed, for $n\ne 0$ (gradient descent is recovered for $n=0$), by simplifying the variable $v$, both schemes can be written in one line:
\vspace{-2mm}
\begin{multline*}
\tag{EE}
    x_{k+1} = x_k + (1-q\sqrt{s})(x_k-x_{k-1}) -m\sqrt{s}\nabla f(x_k) \\
    + ((1-q\sqrt{s})m\sqrt{s}-ns)\nabla f(x_{k-1});
\end{multline*}
\vspace{-7mm}
\begin{multline*}
    \tag{SIE}
x_{k+1} =  x_k + (1-q\sqrt{s}) (x_k - x_{k-1}) \\
- \left(m\sqrt{s}+ns\right) \nabla f(x_k) + (1-q\sqrt{s})m\sqrt{s} \nabla f(x_{k-1}).
\end{multline*}
Remarkably, different choices of parameters yield a rich set of momentum methods, and the reader can probably already notice some configurations which recover well-known optimizers~(see introduction). We explore this in the next subsection.

\subsection{Equivalence between SIE and EE}
\vspace{-2mm}
We show that algorithms obtained from semi-implicit discretization of an accelerated flow can also be seen as explicit discretization of a different accelerated flow.
\begin{tcolorbox}
\begin{lem}[Equivalence between SIE and EE] 
For $n=0$ both EE and SIE reduce to gradient descent. For $n\ne0$, consider parameters $(m_\sie, n_\sie, q)$ and set
\vspace{-6mm}
\begin{align*}
    &m_\ee = m_\sie + \sqrt{s} n_\sie,\\ 
    &n_\ee = (1-q\sqrt{s})n_\sie.
\end{align*}
 EE with stepsize $\sqrt{s}>0$ on \ref{eq:our_model} with parameters $(m_\ee, n_\ee, q)$ leads to the same exact algorithm as the one obtained using SIE with stepsize $\sqrt{s}>0$ on \ref{eq:our_model} with parameters $(m_\sie, n_\sie, q)$.
\label{lemma:equivalence}
\end{lem}
\end{tcolorbox}
\vspace{-3mm}
\begin{proof}
We start from the one-line representation. We get the following conditions for $n\ne 0$:
\vspace{-1mm}
\begin{equation*}
    \begin{cases}
    m_\sie \sqrt{s}+sn_\sie = m_\ee \sqrt{s}\\
    (1-q\sqrt{s})m_\sie\sqrt{s} = (1-q\sqrt{s})m_\ee\sqrt{s} - sn_\ee.
    \end{cases}
\end{equation*}
We substitute the first equation into the second.
\end{proof}
\vspace{-2mm}
As a crucial consequence of the last lemma, Heavy-ball and Nesterov method can be seen both as semi-implicit and explicit integrators on \ref{eq:our_model}. This is illustrated in Tb.~\ref{tb:nag_hb}. Since, as it is well known, \ref{eq: nag} is accelerated, Lemma~\ref{lemma:equivalence} shows that both explicit and semi-implicit Euler integrators can lead to acceleration under well-chosen parameters. In the next subsection, we elaborate more on this finding and recover parameters which lead to acceleration for EE and SIE.
\paragraph{An ODE which gives \ref{eq: nag} under the explicit Euler method.} From Tb.~\ref{tb:nag_hb} and Eq.~\ref{eq:QHM_to_jordan}, we get that
$$\ddot{X} + \big(2\sqrt{\mu} + 2(1-\sqrt{\mu s})\sqrt{s}\nabla^2 f(X)\big)\dot{X} + \nabla f(X) = 0$$
leads to \ref{eq: nag} through EE~(choosing $q=2\sqrt{\mu}$), while
$$\ddot{X} + \big(2\sqrt{\mu} +\sqrt{s}\nabla^2 f(X)\big)\dot{X} + \nabla f(X) = 0$$
recovers \ref{eq: nag} through SIE discretization. These parameter choices lead to acceleration~(see Cor.~\ref{cor:nice_config}). Note that the last equation is equivalent\footnote{The careful reader might notice a factor $1+\sqrt{\mu s}$ in front of the gradient for the ODE in~\citet{shi2019acceleration}. This small difference is only due to the particular definition of semi-implicit integration. If one replaces $q\sqrt{s} v_{k}$ in the RHS of SIE with $q\sqrt{s} v_{k+1}$, then we have complete equivalence.} to \ref{eq: nag-hr}, while the first is not known in the literature. However, Thm.~\ref{thm:continuous} ensures that both ODEs are accelerated. This is enough to show that the sketch in Fig.~\ref{fig:our_picture} is correct.

\begin{center}
\begin{table}
\begin{tabular}{ |c||l|l| } 
 \hline
  & EE discretization & SIE discretization \\ 
\hline\hline
 \ref{eq: hb} & $q = (1-\beta)/\sqrt{s}$ & $q = (1-\beta)/\sqrt{s}$ \\
    & $m =\sqrt{s}$  & $m=0$ \\
    & $n = \beta $ & $n=1$ \\
  \hline
 \ref{eq: nag} & $q = (1-\beta)/\sqrt{s}$ & $q = (1-\beta)/\sqrt{s}$ \\
 &$ m = (1+\beta)\sqrt{s}$ &  $m = \sqrt{s}$ \\
    & $n=\beta^2$  & $n = \beta$ \\
    \hline
\end{tabular}
\caption{\small{\ref{eq: hb} and \ref{eq: nag} with any stepsize $s>0$ and momentum $\beta\in(0,1)$ (see definition in the introduction) can be seen as both EE or SIE numerical integrators.}}
\label{tb:nag_hb}
\end{table}
\end{center}
\vspace{-8mm}
\subsection{Semi-implicit Euler is accelerated}
\begin{figure}
    \centering
 	\includegraphics*[width=\linewidth]{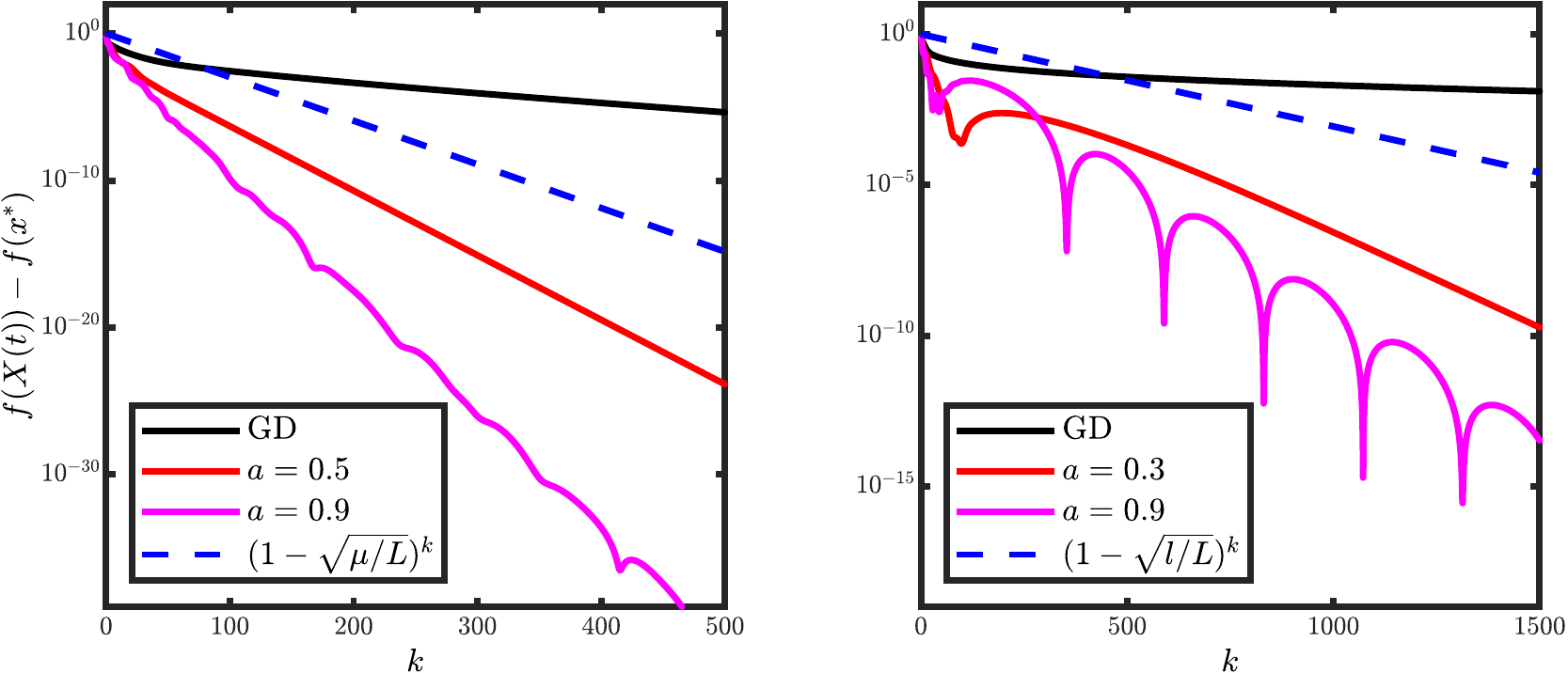} 
 	\caption{\small Convergence of \ref{eq: qhm}. The left plot shows the convergence of $10$-dimensional quadratic with $\mu=0.01$ and $L=1$; the right plot reports $10$-dimensional regularized logistic regression~(random data and labels) with regularization weight $l=10^{-4}$. We specifically used $s=0.5$, $\beta = 1-2\sqrt{\mu s}$ and $s=0.5$, $\beta = 1-2\sqrt{ ls}$, respectively, in the two experiments.}
 	\label{fig: acceleration}
\end{figure}

Leveraging insights from the ODE stability analysis in Thm.~\ref{thm:continuous}
and the lessons learned from semi-implicit Lyapunov function design in recent literature~\citet{shi2018understanding,shi2019acceleration}, our next result establishes a general convergence rate for the semi-implicit Euler method on~\ref{eq:our_model}. In the next subsection, we also provide a similar result for EE, using Lemma~\ref{lemma:equivalence}.
\begin{tcolorbox}
\begin{restatable}[Convergence of SIE]{thm}{SIE}\label{thm:convergence}
	Assume $f$ $L$-smooth and $\mu$-strongly-convex. Let $(x_k)_{k=1}^\infty$ be the sequence obtained from semi-implicit discretization of \ref{eq:our_model} with step $\sqrt{s}$. Let
	\vspace{-2mm}
			\begin{align} \label{eq:parameters_sie}
            0 < m\sqrt{s} \leq \frac{1}{2L},\ 0 < ns\leq m\sqrt{s},\ 0 < q\sqrt{s} \leq \frac{1}{2}.
		\end{align}
	There exists a constant $C>0$ such that, for any $k\in\mathbb{N}$, it holds that 
			\vspace{-2mm}
		\begin{align*}
		&f(x_k) - f(x^*) \leq \left(1+\gamma_2\sqrt{s}\right)^{-k}C, 
		\end{align*}
			\vspace{-2mm}
		where $\gamma_2 :=\frac{1}{5}\min \left(\cfrac{n\mu}{q},\cfrac{q}{1+q^2/(nL)}\right)$.
\end{restatable}
\end{tcolorbox}
\vspace{-2mm}
\begin{proof}[Proof Sketch.] The proof is based on the following energy function inspired by the ODE model ~(cf. Sec.~\ref{sec:cont_model}):
\vspace{-2mm}
\begin{align*}
    \mathcal{E}(k) = &  r_1r_2 (f(x_k) - f(x^*)) -  \frac{r_1r_2m\sqrt{s}}{2} \| \nabla f(x_k)\|^2   \\
    &  + \frac{nr_1^2r_2}{4} \|v_k\|^2 + \frac{1}{4}\|q(x_{k+1}-x^*) - nr_1v_k \|^2, 
\end{align*}
where $r_1 = 1-q\sqrt{s}$, $r_2 = n + mq$ and the last term is a vanishing (as $s\to0$) correction that accounts for the discretization error~(cf.~\citet{shi2019acceleration}). We show $\mathcal{E}(k+1) - \mathcal{E}(k) \leq - \gamma_2\sqrt{s} \mathcal{E}(k+1)$ in App.~\ref{app:SIE}, completing the proof.
\end{proof}
\vspace{-2mm}
The generality of the convergence result allows us to derive accelerated rates for different momentum methods whose convergence rates may even be unknown. We illustrate this by deriving the well-known rate of Nesterov's method in just a few lines. We note that known results on semi-implicit integration such as the ones presented in~\citet{shi2019acceleration} are less general since are limited to high/low resolution or to a fixed viscosity $O(\sqrt{\mu})$.
\vspace{-2mm}
\paragraph{From Thm.~\ref{thm:convergence} to the well-known rate for \ref{eq: nag}.} 
By invoking Thm.~\ref{thm:convergence}, we can recover acceleration of \ref{eq: nag} since it can be written as SIE discretization of~\ref{eq:our_model}~(see Tb.~\ref{tb:nag_hb}). 
\begin{tcolorbox}
\begin{cor}[\ref{eq: nag} is accelerated] 
\label{cor:nice_config} Let $f$ be $L$-smooth and $\mu$-strongly-convex with large\footnotemark{}~condition number $L/\mu\geq 9$. Consider the SIE discretization of \ref{eq:our_model} with $s \leq \frac{1}{4L}$, $q=(1-\beta)/\sqrt{s}$ (with $\beta= 1-2\sqrt{\mu s}$), $m=\sqrt{s}$, $n=\beta$ (i.e. \ref{eq: nag}, see Tb.~\ref{tb:nag_hb}).  The algorithm enjoys the accelerated convergence rate $O((1-\sqrt{\mu/L})^k)$. Namely, $\exists C>0$ such that
\vspace{-2mm}
$$f(x_{k}) - f(x^*) \leq \big(1+\sqrt{\mu s}/15\big)^{-k} C.$$
\end{cor}
\end{tcolorbox}
\footnotetext{The lower bound assumption for conditional number here and in Cor.~\ref{cor:qhm} is purely technical and only serves for a simple illustration of these corollaries.}
\vspace{-2mm}
\begin{proof}
The conditions in in Eq.~\ref{eq:parameters_sie} are satisfied since $s = m\sqrt{s}\le 1/(4L)$, $n=\beta < 1 = \frac{m}{\sqrt{s}}$ and $q\sqrt{s}= 2\sqrt{\mu s} \leq 2 \sqrt{ Ls/9} \leq 1/3$. Thus, $\frac{n\mu}{5q} = \frac{(1-2\sqrt{\mu s})\sqrt{\mu}}{10} \geq \frac{(1-1/3)\sqrt{\mu}}{10}\ge\frac{\sqrt{\mu}}{15}$ and $\frac{q}{5+5q^2/(nL)} \geq \frac{2\sqrt{\mu}}{5+6\mu/L}  \geq \frac{\sqrt{\mu}}{9}$.
\end{proof}
\vspace{-3mm}
\paragraph{From Thm.~4 to a new rate for \ref{eq: qhm}.}
The generality of our model and our discretization analysis provides an accelerated convergence rate for a broad class of momentum methods. Among these methods is quasi-hyperbolic momentum~\citep{ma2018quasi}, which\footnote{\citet{ma2018quasi} presented a normalized second iteration, i.e. $g_{k+1}=bg_k+(1-b)\nabla f(x_k)$, which is generally equivalent to the one we present here by factor rescaling.} shows promises in optimization for neural nets~\citep{choi2019empirical}.
\begin{align}\label{eq: qhm} \tag{QHM}
 	\begin{cases}
 		x_{k+1} = x_k - s((1-a)\nabla f(x_k) + a g_{k+1}) \\
 		g_{k+1} = b g_k +  \nabla f(x_k),
 	\end{cases}
 \end{align}
where $a,b\in(0,1)$. For classification tasks, \ref{eq: qhm} yields an accelerated rate on real-world datasets (even better than \ref{eq: nag}) \citep{ma2018quasi}. Despite empirical benefits, the convergence analysis for this algorithm is limited to quadratics~\citep{gitman2019understanding}. Using Thm.~\ref{thm:convergence}, the next corollary establishes an accelerated rate for \ref{eq: qhm}~(proof in the appendix). 
\begin{tcolorbox}
\begin{restatable}[Convergence of \ref{eq: qhm}]{cor}{QHM} \label{cor:qhm}
Let $f$ be $L$-smooth and $\mu$-strongly-convex with $L/\mu\geq 9$. The iterates of  enjoy a linear convergence rate for $s\leq\frac{1}{4L}$ and $a\leq 1/2$. In particular,  also enjoys convergence rate $O((1-\sqrt{\mu/L})^k)$ for $b = 1 - 2\sqrt{\mu s}$. Namely, $\exists C>0$ such that
\vspace{-2mm}
$$f(x_k) - f(x^*) \leq \Big(1+a\sqrt{\mu s}/10\Big)^{-k} C.$$
\end{restatable} 
\end{tcolorbox}
Fig.~\ref{fig: acceleration} shows the accelerated rate established in the corollary, and its dependency on the parameter $a$. We leave the extension to the stochastic case~(possible with the methodology in~\citet{assran2020convergence}) to future work, for the sake of continuing our discussion on numerical integration.
\vspace{-3mm}
\paragraph{Thm.~\ref{thm:convergence} fails to prove accelerated rate for \ref{eq: hb}.}
An interesting question may arise as a consequence of our results: since \ref{eq: hb} can be recast as semi-implicit discretization of \ref{eq:our_model}, then does invoking Thm.~\ref{thm:convergence} produce a global acceleration proof for \ref{eq: hb}? The answer is \textit{no}, since the convergence result in Thm.~\ref{thm:convergence} is conditioned on $m > 0$; while one needs to set $m=0$ to obtain \ref{eq: hb} by SIE integration. This is not surprising since the Lyapunov functions used in the literature to prove acceleration for \ref{eq: nag} often differ from the one used for convergence of \ref{eq: hb} (see Eq. 3.3 in~\citet{shi2018understanding}). Nonetheless, It is possible to construct an analogue of Thm. 3, using a different Lyapunov function, to derive (non accelerated) convergence for an \ref{eq: hb}-like method.
\vspace{-3mm}
\paragraph{The trade-off speed-stability.} As noted in Remark~\ref{rk:accuracy_stability}, in continuous time one can increase either $m$ or $n$ to infinity and get an arbitrarily fast rate. Thm.~\ref{thm:convergence} shows why a similar phenomenon is not possible in discrete time~(would violate the lower bound in~\citet{nemirovsky1983problem}): for a specific discretization step-size $\sqrt{s}$, Eq.~\ref{eq:parameters_sie} gives us a bound on the maximum $m$ and $n$ we can choose to have guaranteed stability. In other words, if we choose a large value for either $m$ and $n$ to get a faster rate, we would end up with a slow algorithm since numerical stability would require a very small integration step-size. Hence, as expected by the classic theory of convex optimization~\citep{nemirovsky1983problem}, there is a sweet spot which yields $\gamma_2 = O(\sqrt{\mu/L})$ --- a.k.a acceleration. 

\subsection{Explicit Euler is \textit{also} accelerated!}
In the last subsection, we provided a convergence rate for semi-implicit discretization of~\ref{eq:our_model} and showed how this general result can be applied to derive (old and new) convergence rates for momentum methods. However, as already noted a few times, Lemma~\ref{lemma:equivalence} implies that an equivalent theorem can be written for the explicit Euler method.
\allowdisplaybreaks
\begin{tcolorbox}
\begin{restatable}[Convergence of EE]{cor}{EE_from_SIE}\label{thm:convergence_ee}
	Assume $f$ $L$-smooth and $\mu$-strongly-convex. Let $(x_k)_{k=1}^\infty$ be the sequence obtained from semi-implicit discretization of \ref{eq:our_model} with step $\sqrt{s}$. Let
	\vspace{-2mm}
	\begin{align} \label{eq:parameters_ee}
	    \begin{split}
            & 0 < m\sqrt{s} - ns/(1-q\sqrt{s}) \leq \frac{1}{2L}, \\
            & 0 < ns\leq \frac{1-q\sqrt{s}}{2} m\sqrt{s},\quad
             0 < q\sqrt{s} \leq \frac{1}{2}.
        \end{split}
	\end{align}
	There exists a constant $C>0$ such that, for any $k\in\mathbb{N}$, it holds that 
		\begin{align*}
		&f(x_k) - f(x^*) \leq \left(1+\gamma_3\sqrt{s}\right)^{-k}C, 
		\end{align*}
		where $\gamma_3 :=\frac{1}{5}\min \left(\cfrac{n\mu}{q(1-\mu\sqrt{s})},\cfrac{q}{1+q^2/(nL)}\right)$.
\end{restatable}
\end{tcolorbox}
\vspace{-2mm}
\begin{proof}
Consider an explicit method with parameters $(m_\ee, n_\ee, q)$ and a semi-implicit method with parameters $(m_\sie, n_\sie, q)$. Thm.~\ref{thm:convergence} holds if $0<m_\sie\sqrt{s} \leq 1/(2L), 0<s n_\sie \leq  m_\sie\sqrt{s}$ and $q\sqrt{s} \leq 1/2$,
then it is convergent. By Lemma~\ref{lemma:equivalence}, we can recover the parameter of an equivalent explicit method by setting $n_\ee = (1-q\sqrt{s})n_\sie$ and $m_\ee = m_\sie + \sqrt{s} n_\sie$. Combining these conditions with the theorem requirements on $n_\sie$, we get:
$$0<\frac{s n_\ee}{1-q\sqrt{s}} \leq  m_\sie\sqrt{s} = m_\ee\sqrt{s} - \frac{s n_\ee}{1-q\sqrt{s}},$$
which implies the condition on $n_\ee$. For the condition on $m_\ee$, just note that the condition on $m_\sie$ from Thm.~\ref{thm:convergence} implies $\sqrt{s}m_\sie = \sqrt{s}m_\ee -  \frac{s}{1-q\sqrt{s}} n_\ee\le\frac{1}{2L}.$
\vspace{-2mm}
\end{proof}

\begin{figure}[ht]
     \centering
 	\includegraphics*[width=\linewidth]{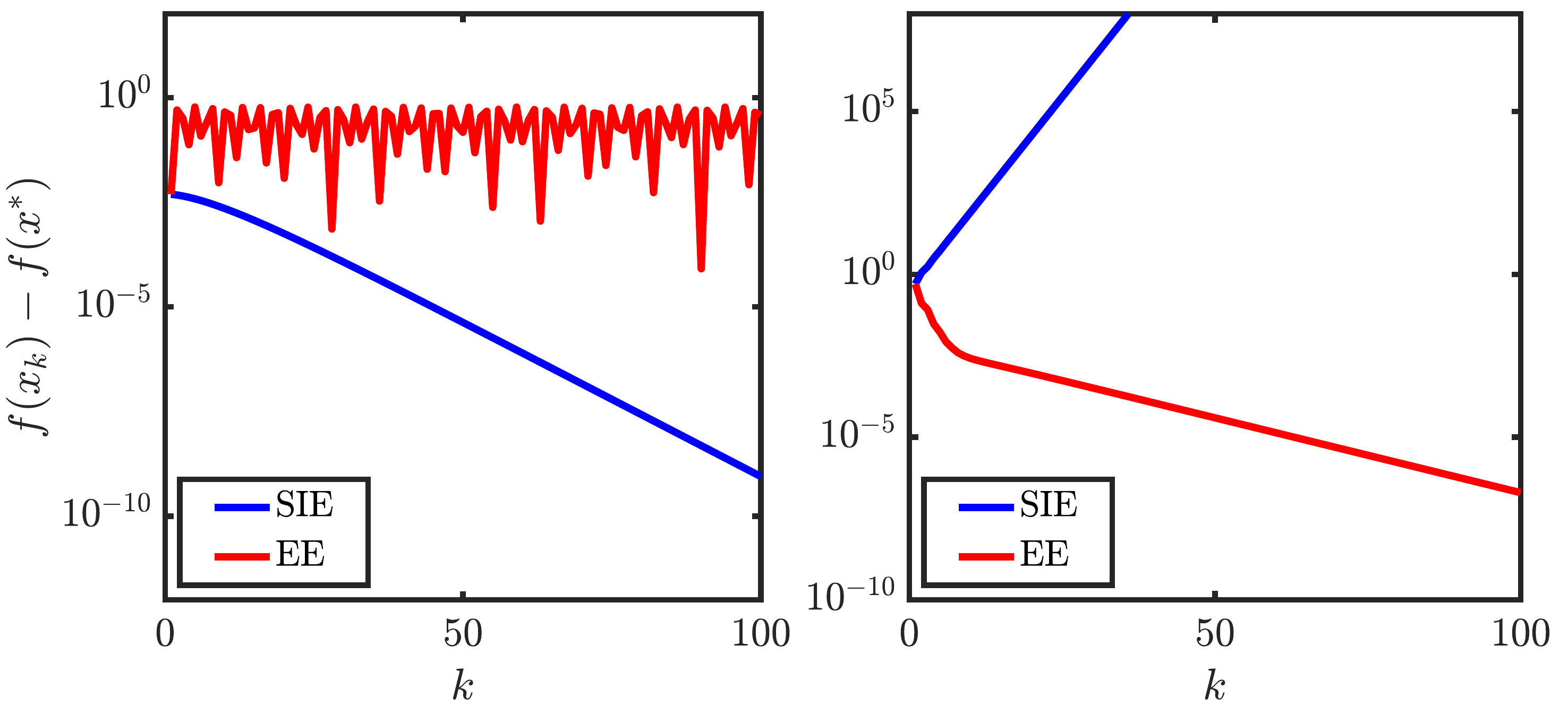} 
 	\caption{\small EE vs. SIE. To show that SIE and EE are neither superior nor inferior to each other, in each subplot, we use the same parameters $m,n,q$ for both SIE and EE discretization. We observe very different behaviours. This suggests the stability and convergence is determined by the joint choice of parameters and numerical integrator together. The objective function here is a 2-dimensional quadratic with $\mu = 0.01, L=1$ and the step-size is $s=1$. In the left plot we use $m = \sqrt{s}, n = 1$ and $q=2\sqrt{\mu}$ and in the right plot we use $m = 2\sqrt{s}, n = 0.5$ and $q=2\sqrt{\mu}$.}
 	\label{fig: contrast}
\end{figure}

\vspace{-2mm}
\paragraph{Stability of EE and SIE. }
For the integration of Hamiltonian systems, semi-implicit Euler is provably more stable than explicit Euler~\citep{hairer2006geometric}. For example, a linearized pendulum integrated with explicit Euler diverges in phase space, while the semi-implicit Euler method is stable and conserves the structure of the ODE system~(energy, volume). In Fig.~\ref{fig: contrast}, we show that for a dissipative (hence not Hamiltonian) system such as \ref{eq:our_model} the situation can be very different: in complete agreement with our equivalence result in Lemma~\ref{lemma:equivalence}, there exists parameter configurations for which EE is stable but SIE is not, and vice versa.

\vspace{-2mm}

\vspace{-2mm}
\section{Behaviour of the discretization error}
\vspace{-2mm}
\label{sec:truncation}
In the last sections, we studied the properties of explicit and semi-implicit integration of \ref{eq:our_model} and showed that both can lead to acceleration. Yet, most recent literature~\citep{wilson2016lyapunov,shi2019acceleration,muehlebach2019dynamical,muehlebach2020optimization} claims that semi-implicit integration \textit{is somehow more natural} for the approximation of partitioned dissipative systems such as \ref{eq:our_model}. Indeed, recent works~\citep{francca2020dissipative,muehlebach2020optimization} showed that the geometric properties of semi-implicit methods combined with backward error analysis~\citep{hairer2006geometric} can be used to successfully prove the preservation of continuous-time rates of convergence up to a controlled error. Instead, our results in Thm.~\ref{thm:convergence_ee} show that explicit Euler discretization --- of a proper ODE --- also leads to an accelerated method~(see also Tb.~\ref{tb:nag_hb}). To conclude our study, we compare semi-implicit and explicit Euler in terms of their approximation error, specifically for the integration of \ref{eq:our_model}. For this particular ODE, EE suffers from a worse local discretization error compared to SIE for the general choice of parameters.
Under particular choices of parameters, EE and SIE yield contractive algorithms. In this case, the error of the both discretization schemes decays exponentially fast. 
\vspace{-3mm}
\paragraph{A trap: local error analysis for the general case.} Consider the following discretization errors:
\begin{align}
    & \Delta_k^{(\ee)} := \| X(k\sqrt{s}) - x^{(\ee)}_{k} \|, \\ 
    & \Delta_k^{(\sie)} := \| X(k\sqrt{s}) - x^{(\sie)}_{k+1} \|.
\end{align}
We compare the above errors for $k=1$ (for one step). Proof/details are provided in the appendix.
\begin{tcolorbox}
\begin{restatable}{lem}{trunError}\label{lemma: trunc_error}
Let $f$ be $L$-smooth and of class $C^2$. If $m=O(\sqrt{s})$, then $\Delta_1^{(\sie)} = O(s^{\sfrac{3}{2}})$ and $\Delta_1^{(\ee)} = O(s)$.
\end{restatable} 
\end{tcolorbox}
\vspace{-1mm}
The above lemma holds for any finite choice of the parameters, and shows that SIE provides a better one-step integration error in the position variable\footnote{It is well known ~\citep{hairer2006geometric} that these methods actually have the same order, since they are $O(s)$ in the velocity.}. This result may lead to a wrong conclusion: semi-implicit integration leads to faster algorithm when discretizing \ref{eq:our_model}. However, this analysis does not provide us a complete picture. Indeed, as we proved in the last section, explicit discretization \textit{can also lead to acceleration} --- in particular, it can recover Nesterov's method. To provide some intuition on why a local error analysis leads to misleading conclusions, we provide a tighter analysis of the integration error for a narrowed set of parameters in \ref{eq:our_model}. 
\vspace{-3mm}
\paragraph{Analysis for contractive cases.} A line of recent works around the connection between acceleration and numerical integration~\citep{orvieto2019shadowing,muehlebach2020optimization,francca2020dissipative} studied the behavior of the discretization error of \ref{eq: nag_ode} as $k\to\infty$, showing interesting \textit{shadowing}\footnote{That is, the discretization bound does \textit{not} explode exponentially due to error accumulation~\citep{chow1994shadowing} if the objective is convex, due to the contraction provided by the landscape.} properties.  The main idea behind shadowing is studying the discretization error when the choice of parameters leads to a contractive algorithm. In this case, one can provide a tighter analysis for the discretization error. The next lemma proves that the integration error of semi-implicit and explicit Euler discretization of \ref{eq:our_model} decays exponentially fast if one properly chooses the parameters.  
\begin{tcolorbox}
\begin{restatable}{lem}{globalError} \label{lem:global_disc_bound}
Let $f$ be $\mu$-strongly-convex and $L$-smooth. For EE discretization of \ref{eq:our_model} obeying Eq.~\ref{eq:parameters_ee}, the discretization error decays as
$
\Delta_k^{(\text{\ee})} = O((1+\gamma_3 \sqrt{s})^{-k})
$ where $\gamma_3$ is defined in Thm.~\ref{thm:convergence_ee}. Furthermore, SIE also enjoys $\Delta_k^{(\text{\sie})} = O((1+\gamma_2 \sqrt{s})^{-k})
$ where $\gamma_2$ is defined in Thm.~\ref{thm:convergence} as long as conditions in Eq.~\ref{eq:parameters_sie} are satisfied. 
\end{restatable}
\end{tcolorbox}
The proof of the last lemma is postponed to the appendix. According to this result, SIE and EE discretization have the same asymptotic integration error properties --- under particular choice of parameters. This similarity is also reflected in the convergence rates. 

\section{Conclusion}

In this paper, we proposed a general ODE model of momentum-based methods for optimizing smooth strongly-convex functions. The generality of our model allows to view different old and new momentum methods as semi-implicit or explicit Euler integrators and to establish novel accelerated convergence rates for both integrators. In particular, our new findings overturn the following old notion: explicit Euler is inferior to semi-implicit~(a.k.a symplectic) Euler because of its unstable nature. Instead, we show that the stability of these integrators is tied to the underlying accelerated ODE. At a deeper level, our methodology provides new challenging insights on the link between accelerated optimization, and numerical integration. 

\paragraph{Acknowledgements.}We are grateful for the enlightening discussions with Christian Lubich and Aurelien Lucchi on the connection between integration accuracy and optimization speed.


\bibliography{bib}

\begin{thebibliography}{54}
\providecommand{\natexlab}[1]{#1}
\providecommand{\url}[1]{\texttt{#1}}
\expandafter\ifx\csname urlstyle\endcsname\relax
  \providecommand{\doi}[1]{doi: #1}\else
  \providecommand{\doi}{doi: \begingroup \urlstyle{rm}\Url}\fi

\bibitem[Ahn(2020)]{ahn2020proximal}
Kwangjun Ahn.
\newblock From proximal point method to {Nesterov}'s acceleration.
\newblock \emph{arXiv:2005.08304}, 2020.

\bibitem[Alecsa(2020)]{alecsa2020long}
Cristian~Daniel Alecsa.
\newblock The long time behavior and the rate of convergence of symplectic
  convex algorithms obtained via splitting discretizations of inertial damping
  systems.
\newblock \emph{arXiv preprint arXiv:2001.10831}, 2020.

\bibitem[Alecsa et~al.(2019)Alecsa, L{\'a}szl{\'o}, and
  Pinta]{alecsa2019extension}
Cristian~Daniel Alecsa, Szil{\'a}rd~Csaba L{\'a}szl{\'o}, and Titus Pinta.
\newblock An extension of the second order dynamical system that models
  {Nesterov’s} convex gradient method.
\newblock \emph{arXiv preprint arXiv:1908.02574}, 2019.

\bibitem[Alimisis et~al.(2020)Alimisis, Orvieto, B{\'e}cigneul, and
  Lucchi]{alimisis2020continuous}
Foivos Alimisis, Antonio Orvieto, Gary B{\'e}cigneul, and Aurelien Lucchi.
\newblock A continuous-time perspective for modeling acceleration in
  {Riemannian} optimization.
\newblock In \emph{International Conference on Artificial Intelligence and
  Statistics}, pages 1297--1307. PMLR, 2020.

\bibitem[Allen-Zhu and Orecchia(2014)]{allen2014linear}
Zeyuan Allen-Zhu and Lorenzo Orecchia.
\newblock Linear coupling: An ultimate unification of gradient and mirror
  descent.
\newblock \emph{arXiv preprint arXiv:1407.1537}, 2014.

\bibitem[Alvarez(2000)]{alvarez2000minimizing}
Felipe Alvarez.
\newblock On the minimizing property of a second order dissipative system in
  {Hilbert} spaces.
\newblock \emph{SIAM Journal on Control and Optimization}, 38\penalty0
  (4):\penalty0 1102--1119, 2000.

\bibitem[Assran and Rabbat(2020)]{assran2020convergence}
Mahmoud Assran and Michael Rabbat.
\newblock On the convergence of {Nesterov's} accelerated gradient method in
  stochastic settings.
\newblock \emph{arXiv preprint arXiv:2002.12414}, 2020.

\bibitem[Attouch and Alvarez(2000)]{attouch2000heavy2}
Hedy Attouch and Felipe Alvarez.
\newblock The heavy ball with friction dynamical system for convex constrained
  minimization problems.
\newblock In \emph{Optimization}, pages 25--35. Springer, 2000.

\bibitem[Attouch et~al.(2000)Attouch, Goudou, and Redont]{attouch2000heavy}
Hedy Attouch, Xavier Goudou, and Patrick Redont.
\newblock The heavy-ball with friction method, i. the continuous dynamical
  system: global exploration of the local minima of a real-valued function by
  asymptotic analysis of a dissipative dynamical system.
\newblock \emph{Communications in Contemporary Mathematics}, 2\penalty0
  (01):\penalty0 1--34, 2000.

\bibitem[B{\'e}gout et~al.(2015)B{\'e}gout, Bolte, and
  Jendoubi]{begout2015damped}
Pascal B{\'e}gout, J{\'e}r{\^o}me Bolte, and Mohamed~Ali Jendoubi.
\newblock On damped second-order gradient systems.
\newblock \emph{Journal of Differential Equations}, 259\penalty0 (7):\penalty0
  3115--3143, 2015.

\bibitem[Benettin and Giorgilli(1994)]{benettin1994hamiltonian}
Giancarlo Benettin and Antonio Giorgilli.
\newblock On the {Hamiltonian} interpolation of near-to-the identity symplectic
  mappings with application to symplectic integration algorithms.
\newblock \emph{Journal of Statistical Physics}, 74\penalty0 (5-6):\penalty0
  1117--1143, 1994.

\bibitem[Betancourt et~al.(2018)Betancourt, Jordan, and
  Wilson]{betancourt2018symplectic}
Michael Betancourt, Michael~I Jordan, and Ashia~C Wilson.
\newblock On symplectic optimization.
\newblock \emph{arXiv preprint arXiv:1802.03653}, 2018.

\bibitem[Bravetti et~al.(2017)Bravetti, Cruz, and Tapias]{bravetti2017contact}
Alessandro Bravetti, Hans Cruz, and Diego Tapias.
\newblock Contact {Hamiltonian} mechanics.
\newblock \emph{Annals of Physics}, 376:\penalty0 17--39, 2017.

\bibitem[Bravetti et~al.(2019)Bravetti, Daza-Torres, Flores-Arguedas, and
  Betancourt]{bravetti2019optimization}
Alessandro Bravetti, Maria~L Daza-Torres, Hugo Flores-Arguedas, and Michael
  Betancourt.
\newblock Optimization algorithms inspired by the geometry of dissipative
  systems.
\newblock \emph{arXiv preprint arXiv:1912.02928}, 2019.

\bibitem[Cabot et~al.(2009)Cabot, Engler, and Gadat]{cabot2009long}
Alexandre Cabot, Hans Engler, and S{\'e}bastien Gadat.
\newblock On the long time behavior of second order differential equations with
  asymptotically small dissipation.
\newblock \emph{Transactions of the American Mathematical Society},
  361\penalty0 (11):\penalty0 5983--6017, 2009.

\bibitem[Choi et~al.(2019)Choi, Shallue, Nado, Lee, Maddison, and
  Dahl]{choi2019empirical}
Dami Choi, Christopher~J Shallue, Zachary Nado, Jaehoon Lee, Chris~J Maddison,
  and George~E Dahl.
\newblock On empirical comparisons of optimizers for deep learning.
\newblock \emph{arXiv preprint arXiv:1910.05446}, 2019.

\bibitem[Chow and Van~Vleck(1994)]{chow1994shadowing}
Shui-Nee Chow and Erik~S Van~Vleck.
\newblock A shadowing lemma approach to global error analysis for initial value
  odes.
\newblock \emph{SIAM Journal on Scientific Computing}, 15\penalty0
  (4):\penalty0 959--976, 1994.

\bibitem[de~Le{\'o}n and Lainz~Valc{\'a}zar(2019)]{de2019contact}
Manuel de~Le{\'o}n and Manuel Lainz~Valc{\'a}zar.
\newblock Contact {Hamiltonian} systems.
\newblock \emph{Journal of Mathematical Physics}, 60\penalty0 (10):\penalty0
  102902, 2019.

\bibitem[Defazio(2019)]{defazio2019curved}
Aaron Defazio.
\newblock On the curved geometry of accelerated optimization.
\newblock In \emph{Advances in Neural Information Processing Systems}, pages
  1764--1773, 2019.

\bibitem[Diakonikolas and Jordan(2019)]{diakonikolas2019generalized}
Jelena Diakonikolas and Michael~I Jordan.
\newblock Generalized momentum-based methods: A {Hamiltonian} perspective.
\newblock \emph{arXiv preprint arXiv:1906.00436}, 2019.

\bibitem[Diakonikolas and Orecchia(2017)]{diakonikolas2017accelerated}
Jelena Diakonikolas and Lorenzo Orecchia.
\newblock Accelerated extra-gradient descent: A novel accelerated first-order
  method.
\newblock \emph{arXiv preprint arXiv:1706.04680}, 2017.

\bibitem[Fazlyab et~al.(2018)Fazlyab, Ribeiro, Morari, and
  Preciado]{fazlyab2018analysis}
Mahyar Fazlyab, Alejandro Ribeiro, Manfred Morari, and Victor~M Preciado.
\newblock Analysis of optimization algorithms via integral quadratic
  constraints: Nonstrongly convex problems.
\newblock \emph{SIAM Journal on Optimization}, 28\penalty0 (3):\penalty0
  2654--2689, 2018.

\bibitem[Flammarion and Bach(2015)]{flammarion2015averaging}
Nicolas Flammarion and Francis Bach.
\newblock From averaging to acceleration, there is only a step-size.
\newblock In \emph{Conference on Learning Theory}, pages 658--695, 2015.

\bibitem[Fran{\c{c}}a et~al.(2020{\natexlab{a}})Fran{\c{c}}a, Jordan, and
  Vidal]{francca2020dissipative}
Guilherme Fran{\c{c}}a, Michael~I Jordan, and Ren{\'e} Vidal.
\newblock On dissipative symplectic integration with applications to
  gradient-based optimization.
\newblock \emph{arXiv preprint arXiv:2004.06840}, 2020{\natexlab{a}}.

\bibitem[Fran{\c{c}}a et~al.(2020{\natexlab{b}})Fran{\c{c}}a, Sulam, Robinson,
  and Vidal]{francca2020conformal}
Guilherme Fran{\c{c}}a, Jeremias Sulam, Daniel Robinson, and Ren{\'e} Vidal.
\newblock Conformal symplectic and relativistic optimization.
\newblock \emph{Advances in Neural Information Processing Systems}, 33,
  2020{\natexlab{b}}.

\bibitem[Gavurin(1958)]{gavurin1958nonlinear}
Mark~Konstantinovich Gavurin.
\newblock Nonlinear functional equations and continuous analogues of iteration
  methods.
\newblock \emph{Izvestiya Vysshikh Uchebnykh Zavedenii. Matematika}, pages
  18--31, 1958.

\bibitem[Gitman et~al.(2019)Gitman, Lang, Zhang, and
  Xiao]{gitman2019understanding}
Igor Gitman, Hunter Lang, Pengchuan Zhang, and Lin Xiao.
\newblock Understanding the role of momentum in stochastic gradient methods.
\newblock In \emph{Advances in Neural Information Processing Systems}, pages
  9630--9640, 2019.

\bibitem[Hairer(1994)]{hairer1994backward}
Ernst Hairer.
\newblock Backward analysis of numerical integrators and symplectic methods.
\newblock \emph{Annals of Numerical Mathematics}, 1:\penalty0 107--132, 1994.

\bibitem[Hairer et~al.(2006)Hairer, Lubich, and Wanner]{hairer2006geometric}
Ernst Hairer, Christian Lubich, and Gerhard Wanner.
\newblock \emph{Geometric numerical integration: structure-preserving
  algorithms for ordinary differential equations}, volume~31.
\newblock Springer Science \& Business Media, 2006.

\bibitem[Hu and Lessard(2017)]{hu2017dissipativity}
Bin Hu and Laurent Lessard.
\newblock Dissipativity theory for {Nesterov}'s accelerated method.
\newblock In \emph{Proceedings of the 34th International Conference on Machine
  Learning-Volume 70}, pages 1549--1557. JMLR. org, 2017.

\bibitem[Khalil and Grizzle(2002)]{khalil2002nonlinear}
Hassan~K Khalil and Jessy~W Grizzle.
\newblock \emph{Nonlinear systems}, volume~3.
\newblock Prentice Hall Upper Saddle River, NJ, 2002.

\bibitem[Krichene et~al.(2015)Krichene, Bayen, and
  Bartlett]{krichene2015accelerated}
Walid Krichene, Alexandre Bayen, and Peter~L Bartlett.
\newblock Accelerated mirror descent in continuous and discrete time.
\newblock \emph{Advances in neural information processing systems},
  28:\penalty0 2845--2853, 2015.

\bibitem[Lessard et~al.(2016)Lessard, Recht, and Packard]{lessard2016analysis}
Laurent Lessard, Benjamin Recht, and Andrew Packard.
\newblock Analysis and design of optimization algorithms via integral quadratic
  constraints.
\newblock \emph{SIAM Journal on Optimization}, 26\penalty0 (1):\penalty0
  57--95, 2016.

\bibitem[Lubich(2008)]{lubich2008quantum}
Christian Lubich.
\newblock \emph{From quantum to classical molecular dynamics: reduced models
  and numerical analysis}.
\newblock European Mathematical Society, 2008.

\bibitem[Ma and Yarats(2018)]{ma2018quasi}
Jerry Ma and Denis Yarats.
\newblock Quasi-hyperbolic momentum and adam for deep learning.
\newblock \emph{arXiv preprint arXiv:1810.06801}, 2018.

\bibitem[McLachlan and Perlmutter(2001)]{mclachlan2001conformal}
Robert McLachlan and Matthew Perlmutter.
\newblock Conformal {Hamiltonian} systems.
\newblock \emph{Journal of Geometry and Physics}, 39\penalty0 (4):\penalty0
  276--300, 2001.

\bibitem[McLachlan and Quispel(2002)]{mclachlan2002splitting}
Robert~I McLachlan and G~Reinout~W Quispel.
\newblock Splitting methods.
\newblock \emph{Acta Numerica}, 11:\penalty0 341, 2002.

\bibitem[Muehlebach and Jordan(2019)]{muehlebach2019dynamical}
Michael Muehlebach and Michael~I Jordan.
\newblock A dynamical systems perspective on {Nesterov} acceleration.
\newblock \emph{arXiv preprint arXiv:1905.07436}, 2019.

\bibitem[Muehlebach and Jordan(2020)]{muehlebach2020optimization}
Michael Muehlebach and Michael~I Jordan.
\newblock Optimization with momentum: Dynamical, control-theoretic, and
  symplectic perspectives.
\newblock \emph{arXiv preprint arXiv:2002.12493}, 2020.

\bibitem[Nemirovsky and Yudin(1983)]{nemirovsky1983problem}
Arkadi~Semenovich Nemirovsky and David~Borisovich Yudin.
\newblock \emph{Problem complexity and method efficiency in optimization.}
\newblock Wiley, 1983.

\bibitem[Nesterov(1983)]{nesterov1983method}
Yurii~E Nesterov.
\newblock A method for solving the convex programming problem with convergence
  rate {$O(1/k^2)$}.
\newblock In \emph{Dokl. akad. nauk Sssr}, volume 269, pages 543--547, 1983.

\bibitem[Orvieto and Lucchi(2019)]{orvieto2019shadowing}
Antonio Orvieto and Aurelien Lucchi.
\newblock Shadowing properties of optimization algorithms.
\newblock In \emph{Advances in Neural Information Processing Systems}, pages
  12692--12703, 2019.

\bibitem[Orvieto et~al.(2020)Orvieto, Kohler, and Lucchi]{orvieto2020role}
Antonio Orvieto, Jonas Kohler, and Aurelien Lucchi.
\newblock The role of memory in stochastic optimization.
\newblock In \emph{Uncertainty in Artificial Intelligence}, pages 356--366.
  PMLR, 2020.

\bibitem[Polyak(1964)]{polyak1964some}
Boris~T Polyak.
\newblock Some methods of speeding up the convergence of iteration methods.
\newblock \emph{USSR Computational Mathematics and Mathematical Physics},
  4\penalty0 (5):\penalty0 1--17, 1964.

\bibitem[Sanz-Serna and Zygalakis(2020)]{sanz2020connections}
JM~Sanz-Serna and Konstantinos~C Zygalakis.
\newblock The connections between {Lyapunov} functions for some optimization
  algorithms and differential equations.
\newblock \emph{arXiv preprint arXiv:2009.00673}, 2020.

\bibitem[Sanz~Serna and Zygalakis(2020)]{sanz2020contractivity}
JM~Sanz~Serna and Konstantinos~C Zygalakis.
\newblock Contractivity of {Runge}--{Kutta} methods for convex gradient
  systems.
\newblock \emph{SIAM Journal on Numerical Analysis}, 58\penalty0 (4):\penalty0
  2079--2092, 2020.

\bibitem[Shi et~al.(2018)Shi, Du, Jordan, and Su]{shi2018understanding}
Bin Shi, Simon~S Du, Michael~I Jordan, and Weijie~J Su.
\newblock Understanding the acceleration phenomenon via high-resolution
  differential equations.
\newblock \emph{arXiv preprint arXiv:1810.08907}, 2018.

\bibitem[Shi et~al.(2019)Shi, Du, Su, and Jordan]{shi2019acceleration}
Bin Shi, Simon~S Du, Weijie Su, and Michael~I Jordan.
\newblock Acceleration via symplectic discretization of high-resolution
  differential equations.
\newblock In \emph{Advances in Neural Information Processing Systems}, pages
  5744--5752, 2019.

\bibitem[Su et~al.(2014)Su, Boyd, and Candes]{su2014differential}
Weijie Su, Stephen Boyd, and Emmanuel Candes.
\newblock A differential equation for modeling {Nesterov}’s accelerated
  gradient method: Theory and insights.
\newblock In \emph{Advances in Neural Information Processing Systems}, pages
  2510--2518, 2014.

\bibitem[Wibisono et~al.(2016)Wibisono, Wilson, and
  Jordan]{wibisono2016variational}
Andre Wibisono, Ashia~C Wilson, and Michael~I Jordan.
\newblock A variational perspective on accelerated methods in optimization.
\newblock \emph{Proceedings of the National Academy of Sciences}, 113\penalty0
  (47):\penalty0 E7351--E7358, 2016.

\bibitem[Wilson et~al.(2016)Wilson, Recht, and Jordan]{wilson2016lyapunov}
Ashia~C Wilson, Benjamin Recht, and Michael~I Jordan.
\newblock A {Lyapunov} analysis of momentum methods in optimization.
\newblock \emph{arXiv preprint arXiv:1611.02635}, 2016.

\bibitem[Xu et~al.(2018)Xu, Wang, and Gu]{xu2018continuous}
Pan Xu, Tianhao Wang, and Quanquan Gu.
\newblock Continuous and discrete-time accelerated stochastic mirror descent
  for strongly convex functions.
\newblock In \emph{International Conference on Machine Learning}, pages
  5492--5501, 2018.

\bibitem[Zhang et~al.(2018)Zhang, Mokhtari, Sra, and
  Jadbabaie]{zhang2018direct}
Jingzhao Zhang, Aryan Mokhtari, Suvrit Sra, and Ali Jadbabaie.
\newblock Direct {Runge}-{Kutta} discretization achieves acceleration.
\newblock In \emph{Advances in neural information processing systems}, pages
  3900--3909, 2018.

\bibitem[Zhang et~al.(2019)Zhang, Sra, and Jadbabaie]{zhang2019acceleration}
Jingzhao Zhang, Suvrit Sra, and Ali Jadbabaie.
\newblock Acceleration in first order quasi-strongly convex optimization by ode
  discretization.
\newblock In \emph{2019 IEEE 58th Conference on Decision and Control (CDC)},
  pages 1501--1506. IEEE, 2019.

\end{thebibliography}
\bibliographystyle{plainnat}

\onecolumn
\newpage
\newpage
\renewcommand\thesection{\Alph{section}}
\setcounter{section}{0}
\begin{center}
	\textbf{\Large{Appendix: Proofs and Supplementaries}}
\end{center}

\section{Proof for Theorem~\ref{thm:continuous}}
For convenience of the reader, we report here our generalized model for momentum methods \eqref{eq:our_model}, motivated in the main paper.
\begin{align*}
\begin{cases}
\dot{X} = - m \nabla f(X) - nV \\
\dot{V} = \nabla f(X) - qV.  \tag{GM-ODE}
\end{cases}
\end{align*}
\begin{tcolorbox}
\continuous*
\end{tcolorbox}

\begin{proof}
	We propose the Lyapunov function
	\begin{align} 
		\mathcal{E}(t) = & \underbrace{(qm+n)}_{c_1}\big(f(X(t))-f(x^*)\big) + \underbrace{\frac{n(qm+n)}{4}}_{c_2}\|V(t)\|^2 + \underbrace{\frac{1}{4}}_{c_3}\|q(X(t)-x^*)-nV(t)\|^2,
	\end{align}
	consisting of quadratic and mixing parts
	\begin{align}
	\mathcal{E}_1(t) = f(X(t)) - f(x^*), \quad \mathcal{E}_2(t)  = \|V(t)\|^2, \quad	\mathcal{E}_3(t) = \|-nV(t) + q(X(t)-x^*)\|^2.
	\end{align}
	The derivatives of each quadratic part are
	\begin{align}
	\frac{d}{dt} \mathcal{E}_1(t) & = - m \|\nabla f(X(t)\|^2 - n\langle \nabla f(X(t)), V(t)\rangle
	\end{align}
	and
	\begin{align}
	\frac{d}{dt} \mathcal{E}_2(t) & = - 2q \|V(t)\|^2 + 2\langle \nabla f(X(t)), V(t)\rangle,
	\end{align}
	along with that of the mixing term:
	\begin{align}
	\frac{d}{dt} \mathcal{E}_3(t) = & 2\langle -n\dot{V}(t) + q\dot{X}(t), - n V(t) + q (X(t)-x^*) \rangle \nonumber\\
	= & - 2(qm + n) \langle \nabla f(X(t)), - n V(t) + q (X(t)-x^*) \rangle \nonumber\\
	= & - 2q(qm + n)  \langle \nabla f(X(t)), X(t)-x^* \rangle \nonumber + 2n (qm + n)  \langle \nabla f(X(t)), V(t) \rangle   \nonumber\\
	\leq & - 2q(qm + n)  \Big(f(X(t)) - f(x^*)\Big)- \mu q(qm + n) \|X(t)-x^*\|^2 \nonumber\\ 
	& + 2n(qm + n) \langle \nabla f(X(t)), V(t) \rangle,
	\end{align}
	where last inequality is due to the strong convexity. Plugging the value of $c_1$, $c_2$ and $c_3$, we have
	\begin{align}
	\frac{d}{dt}\mathcal{E}(t) \leq & -\frac{q(n+qm)}{2}\Big(\big(f(X(t)) - f(x^*) \big)  + \frac{\mu}{2}\|X(t)-x^*\|^2 + n \|V(t)\|^2\Big). 
	\end{align}
	Besides, the mixing term can be upper-bounded by  
	\begin{align}
	\mathcal{E}_3(t)  \leq 2q^2 \|X(t)-x^*\|^2 + 2n^2\|V(t)\|^2.
	\end{align}
	Therefore we have $\mathcal{E}(t)$ satisfying 
	\begin{align}
	\mathcal{E}(t) \leq (qm + n) \big(f(X(t))-f(x^*)\big) + q^2\|X(t)-x^*\|^2/2 + \Big(n^2/2+\frac{n(n + qm)}{4}\Big)\|V(t)\|^2, 
	\end{align}
	which implies
	\begin{align}
	\frac{d}{dt} \mathcal{E}(t) \leq - \min\left\{\frac{\mu (n+qm)}{2q}, \frac{q}{2}\right\} \cdot \mathcal{E}(t).
	\end{align}
	We then conclude using Gronwall's lemma~\citep{khalil2002nonlinear}.
\end{proof}

\section{Proof for Theorem~\ref{thm:convergence}}
\label{app:SIE}
For convenience of the reader, we repeat here the semi-implicit integrator of \ref{eq:our_model} we seek to study:
\begin{align*}
\textrm{(SIE)}: \quad\begin{cases}
x_{k+1} - x_k = - m\sqrt{s} \nabla f(x_k) - n \sqrt{s} v_k \\
v_{k+1} - v_k = \sqrt{s} \nabla f(x_{k+1}) - q\sqrt{s} v_{k}.
\end{cases}
\end{align*}
In compact notation, the second iteration can be written as
\begin{align} \label{eq: compact_sie1}
    r_1(v_{k+1} - v_k) = \sqrt{s} \nabla f(x_{k+1}) - q\sqrt{s} v_{k+1}
\end{align}
or
\begin{align} \label{eq: compact_sie2}
    r_1v_k = v_{k+1}   - \sqrt{s} \nabla f(x_{k+1}), 
\end{align}
where $r_1=1-q\sqrt{s}$.

\begin{tcolorbox}
\SIE*
\end{tcolorbox}

\begin{proof}
We propose the discrete Lyapunov function defined as
\begin{align}
    \mathcal{E}(k) = \textcolor{blue}{r_1r_2 (f(x_k) - f(x^*))} &  + \textcolor{red}{\frac{1}{4}\|q(x_{k+1}-x^*) - nr_1v_k \|^2} + \textcolor{magenta}{\frac{nr_1^2r_2}{4} \|v_k\|^2} -  \frac{r_1r_2m\sqrt{s}}{2} \| \nabla f(x_k)\|^2.
\end{align}
We use colors for different parts to keep track of related terms in the derivation. As the first step, thanks to $L$-Lipshitz smoothness, we have
\begin{align} \label{eq:Lsmoothness_result}
    f(x_{k+1}) - f(x_k)  \leq & \langle \nabla f(x_{k+1}), x_{k+1} - x_k \rangle - \frac{1}{2L}\|\nabla f(x_{k+1}) - \nabla f(x_k)\|^2 \nonumber\\
    = & - m\sqrt{s}\langle \nabla f(x_k), \nabla f(x_{k+1}) \rangle - n\sqrt{s} \langle v_k, \nabla f(x_{k+1})\rangle \nonumber\\ 
    & - \frac{1}{2L}\|\nabla f(x_{k+1}) - \nabla f(x_k)\|^2.
\end{align}
We proceed by computing the difference in $\mathcal{E}$ in two subsequent iterations. Denote $r_2 = n + mq$, we have
\allowdisplaybreaks
\begin{align}
\mathcal{E}(k+1) - \mathcal{E}(k) 
    \stackrel{(A)}{\leq} & - \textcolor{blue}{r_1r_2m\sqrt{s}\langle \nabla f(x_k), \nabla f(x_{k+1}) \rangle - r_1r_2n\sqrt{s} \langle v_k, \nabla f(x_{k+1})\rangle} \nonumber\\ 
    & - \textcolor{blue}{\frac{r_1r_2}{2L}\|\nabla f(x_{k+1}) - \nabla f(x_k)\|^2} + \textcolor{red}{ \frac{1}{4} \| q(x_{k+2} - x_{k+1}) - nr_1(v_{k+1} - v_k) \|^2} \nonumber\\
    & + \textcolor{red}{\frac{1}{2} \langle q(x_{k+2} - x_{k+1}) - nr_1 (v_{k+1} - v_k), q(x_{k+1}-x^*) - nv_{k+1}  + n\sqrt{s}\nabla f(x_{k+1}) \rangle} \nonumber\\
    & + \textcolor{magenta}{\frac{nr_1^2r_2}{4}\|v_{k+1}\|^2 - \frac{nr_2}{4}\|v_{k+1} - \sqrt{s}\nabla f(x_{k+1})\|^2} \nonumber\\ 
    & -  \frac{r_1r_2m\sqrt{s}}{2} \Big( \|\nabla f(x_{k+1})\|^2 - \|\nabla f(x_k)\|^2 \Big) \nonumber\\
    \stackrel{(B)}{=} &  \textcolor{blue}{-r_1r_2m\sqrt{s}\langle \nabla f(x_k), \nabla f(x_{k+1}) \rangle - r_1r_2n\sqrt{s} \langle v_k, \nabla f(x_{k+1})\rangle} \nonumber\\ 
    & \textcolor{blue}{- \frac{r_1r_2}{2L}\|\nabla f(x_{k+1}) - \nabla f(x_k)\|^2} - \textcolor{red}{\frac{r_2(2n-r_2)}{4}s \| \nabla f(x_{k+1}) \|^2}\nonumber\\
    & - \textcolor{red}{\frac{r_2}{2} \sqrt{s} \langle  \nabla f(x_{k+1}), q(x_{k+1}-x^*) - nv_{k+1} \rangle} \nonumber\\
    & -\textcolor{magenta}{\frac{nr_2(1-r_1^2)}{4} \|v_{k+1}\|^2 - \frac{nr_2}{4}s\|\nabla f(x_{k+1})\|^2 + \frac{nr_2}{2}\sqrt{s} \langle \nabla f(x_{k+1}), v_{k+1} \rangle} \nonumber\\
    & - \frac{r_1r_2m\sqrt{s}}{2} \Big( \|\nabla f(x_{k+1})\|^2 - \|\nabla f(x_k)\|^2 \Big) \nonumber\\
    \stackrel{(C)}{=} & nr_2\sqrt{s}\langle  \nabla f(x_{k+1}), \textcolor{red}{v_{k+1}/2}+\textcolor{magenta}{v_{k+1}/2} - \textcolor{blue}{r_1v_k} \rangle \nonumber\\
    & +\frac{r_1r_2}{2}m\sqrt{s}\Big(\|\nabla f(x_{k+1})\|^2 -\textcolor{blue}{2\langle \nabla f(x_{k+1}), \nabla f(x_{k}) \rangle } + \|\nabla f(x_k)\|^2\Big) \nonumber\\
    & - \Big(\textcolor{red}{\frac{r_2(2n-r_2)}{4}s} + \textcolor{magenta}{\frac{nr_2}{4}s} + r_1r_2m\sqrt{s} \Big)\| \nabla f(x_{k+1}) \|^2  -\textcolor{magenta}{\frac{nr_2(1-r_1^2)}{4} \|v_{k+1}\|^2} \nonumber \\
    &- \textcolor{blue}{\frac{r_1r_2}{2L}\|\nabla f(x_{k+1}) - \nabla f(x_k)\|^2} - \textcolor{red}{\frac{r_2}{2}q\sqrt{s} \langle  \nabla f(x_{k+1}), x_{k+1} - x^* \rangle}. 
\end{align}
In step (A), we use smoothness of $f$ as stated in Eq.~\ref{eq:Lsmoothness_result} for the blue term. Also, we used the inequality $\| a\|^2 - \| b \|^2 = \| a-b \|^2 + 2 \langle a-b, b \rangle$ where $a = q(x_{k+2} - x^*) -nr_1 v_k$ and $b = q(x_{k+1} -x^*) -nr_1 v_k$ to obtain the red term. In particular,
\begin{align}
a-b &= q(x_{k+2} - x_{k+1})-nr_1(v_{k+1}-v_k) \nonumber\\
&= -mq\sqrt{s}\nabla f(x_{k+1})-nq\sqrt{s}v_{k+1}-n\sqrt{s}\nabla f(x_{k+1}) + nq\sqrt{s} v_{k+1} \nonumber\\
&= -r_2\sqrt{s}\nabla f(x_{k+1}).
\end{align}
In step (B), we incorporate the recurrence of SIE. Step (C) is a simple re-arrangement of terms.

We can easily verify the following identities:
\begin{align}
    \sqrt{s}\langle \nabla f(x_{k+1}),  v_{k+1} - r_1v_k \rangle = s \| \nabla f(x_{k+1})\|^2
\end{align}
and
\begin{align}
    \|\nabla f(x_{k+1})\|^2 -2\langle \nabla f(x_{k+1}), \nabla f(x_{k}) \rangle + \|\nabla f(x_k)\|^2 = \|\nabla f(x_{k+1}) - \nabla f(x_k)\|^2.
\end{align}
We have
\begin{align} \label{eq:layp_mid_bound1}
\mathcal{E}(k+1) - \mathcal{E}(k) \leq &  nr_2s\| \nabla f(x_{k+1}) \|^2 +\frac{r_1r_2}{2}m\sqrt{s} \|\nabla f(x_{k+1}) - \nabla f(x_k)\|^2  \nonumber\\
    & - r_2s\Big(\frac{2n-r_2}{4} + \frac{n}{4} + \frac{r_1m}{\sqrt{s}} \Big)\| \nabla f(x_{k+1}) \|^2  -\frac{nr_2(1-r_1^2)}{4} \|v_{k+1}\|^2 \nonumber\\
    & - \frac{r_1r_2}{2L}\|\nabla f(x_{k+1}) - \nabla f(x_k)\|^2 - \frac{r_2}{2}q\sqrt{s} \langle  \nabla f(x_{k+1}), x_{k+1} - x^* \rangle. 
\end{align}
We leverage $\mu$-strong convexity of $f$ to get 
\begin{align}
    \langle \nabla f(x_{k+1}), x_{k+1} - x^* \rangle \geq f(x_{k+1}) - f(x^*) + \frac{\mu}{2} \| x_{k+1} - x^* \|^2 .
\end{align}
Applying the above inequality to the last term of Eq.~\ref{eq:layp_mid_bound1}, we obtain
\begin{align}
\mathcal{E}(k+1) - \mathcal{E}(k) & \leq  - \frac{r_2}{2}q\sqrt{s} (f(x_{k+1})-f(x^*)) - \frac{r_2\mu}{4}q\sqrt{s}\|x_{k+1}-x^*\|^2 \nonumber\\
    & -\frac{r_1r_2}{2}( 1/L - m\sqrt{s} ) \|\nabla f(x_{k+1}) - \nabla f(x_k)\|^2  - \frac{nr_2(1-r_1^2)}{4} \|v_{k+1}\|^2 \nonumber\\
    & - r_2s\Big(\frac{2n-r_2}{4} + \frac{n}{4} + \frac{r_1m}{\sqrt{s}} - n \Big)\|\nabla f(x_{k+1})\|^2.
\end{align}

Now we plug in the the value of $r_1$, $r_2$ and calculate
\begin{align}
    1 - r_1^2 = 1 - (1 - q\sqrt{s}) ^ 2 = q\sqrt{s}(2 - q\sqrt{s}) \geq q\sqrt{s},  
\end{align}
where we used the condition $q\sqrt{s} \leq 1/2$. Next, since $m\sqrt{s} \leq 1/(2L)$, $n \leq m/\sqrt{s}$ and $r_1 = 1 - q\sqrt{s} \geq 1/2$, it holds that
\begin{align}
    \frac{2n-r_2}{4} + \frac{r_1m}{\sqrt{s}} - \frac{3n}{4} = \frac{n-mq}{4} + \frac{r_1m}{\sqrt{s}} - \frac{3n}{4} = \frac{r_1m}{\sqrt{s}}-\frac{n}{2} - \frac{mq}{4} \geq  - \frac{mq}{4}.
\end{align}
Hence, the difference between two iterations can be upper-bounded as follows:
\begin{align}
    \mathcal{E}(k+1) - \mathcal{E}(k) \leq &
    -\frac{r_2q\sqrt{s}}{2}\Big( f(x_{k+1}) - f(x^*) + \frac{\mu}{2}\|x_{k+1} - x^*\|^2 + n\|v_{k+1}\|^2/2 - \frac{m\sqrt{s}}{2} \|\nabla f(x_{k+1})\|^2 \Big) \nonumber\\
    = & - \frac{r_2q\sqrt{s}}{2}\Big( (1-r_3)[f(x_{k+1}) - f(x^*)] + \frac{\mu}{2}\|x_{k+1} - x^*\|^2   \nonumber\\
    & \qquad + n\|v_{k+1}\|^2/2 + r_3[ f(x_{k+1}) -f(x^*) -\frac{1}{2L}\|\nabla f(x_{k+1})\|^2] \Big),
\end{align}
where $r_3 = Lm\sqrt{s} \leq 1 /2$ and the bound remains legal since $1 - r_3 \geq 1/2$.

On the other hand, our candidate Lyapunov function at iteration $k$ itself can be upper-bounded as 
\begin{align}
    \mathcal{E}(k)  = &\ r_1r_2 (f(x_k)-f(x^*)   + \frac{1}{4}\|q(x_{k+1}-x^*) - nr_1v_k \|^2  + \frac{nr_1^2r_2}{4}\|v_k\|^2 -  \frac{r_1r_2m\sqrt{s}}{2} \| \nabla f(x_k)\|^2 \nonumber\\
    \stackrel{(A)}{=} & \ r_1r_2 (f(x_k)-f(x^*))   + \frac{1}{4}\|q(x_{k}-x^*) - nv_k - mq\sqrt{s}\nabla f(x_k) \|^2 + \frac{nr_1^2r_2}{4}\|v_k\|^2 \nonumber\\
    &  -  r_1r_2m\sqrt{s} \| \nabla f(x_k)\|^2 / 2  \nonumber\\
    \stackrel{(B)}{\leq} & \ r_1r_2 (f(x_k)-f(x^*))   + q^2\|x_k-x^*\|^2 + n^2\|v_k\|^2  +\frac{q^2m^2s}{2}\|\nabla f(x_k) \|^2+\frac{nr_1^2r_2}{4} \|v_k\|^2 \nonumber\\
    & -  r_1r_2m\sqrt{s} \| \nabla f(x_k)\|^2 /2  \nonumber\\
    = & \ r_1r_2(1-r_3 +r_4) (f(x_k)-f(x^*)) + q^2\|x_k-x^*\|^2 + (n^2+nr_1^2r_2/4)\|v_k\|^2  \nonumber\\
    & + r_1r_2(r_3 - r_4)[  f(x_k) -f(x^*)- \frac{1}{2L}\| \nabla f(x_k)\|^2], 
\end{align}
with $r_4 = Lq^2m^2s/(r_1r_2)$. Precisely, step (A) is obtained by replacing SIE update for the term $x_{k+1}$. (B) is obtained by repeatedly using the inequality $\|a+b\|^2\leq 2\|a\|^2+2\|b\|^2$. Finally, noting that $f(x_k) - f(x^*) \geq \frac{1}{2L}\|\nabla f(x_k)\|^2$, we have
\begin{multline}
        \mathcal{E}(k) \leq  r_2\Big(r_1(1-r_3 + r_4) [f(x_k) -f(x^*)] + \frac{q^2}{n}\|x_k - x^*\|^2  + 5n\|v_k\|^2/4 \\
        + r_1(r_3-r_4)[f(x_k)-f(x^*)-\frac{1}{2L}\|\nabla f(x_k)\|^2]\Big),
\end{multline}
since $r_2 = n+mq \geq n$. It is reckoned that $\mathcal{E}(k+1) - \mathcal{E}(k)$ and $\mathcal{E}(k)$ share identical parts except for different coefficients. Now we aim at obtaining following inequality
\begin{align}
    \mathcal{E}(k+1) - \mathcal{E}(k) \leq - \gamma_2\sqrt{s} \mathcal{E}(k+1).
\end{align}
To achieve this, $\gamma_2$ should be the minimal ratio for coefficients of each parts of $\mathcal{E}(k+1) - \mathcal{E}(k)$ to those of $\mathcal{E}(k)$. It is easy then to notice that $\gamma_2$ should be smaller than $q/5$ and $n\mu/(4q)$. Besides it should also hold that
\begin{align}
    \frac{r_2q}{2r_1r_2}\frac{r_3}{r_3 - r_4} \geq \frac{q}{2}\frac{1-r_3}{1-(r_3 - r_4)} = \frac{q}{2}\frac{1-r_3}{1-r_3(1 - r_4/r_3)} \geq \frac{q}{2}\frac{1-1/2}{1-1/2(1 - \frac{q^2}{nL})} \geq \frac{q}{2}\frac{1}{1+ \frac{q^2}{nL}} \geq\gamma_2,
\end{align}
due to the fact $\frac{r_4}{r_3} = \frac{q^2m\sqrt{s}}{r_1r_2} \leq \frac{q^2}{nL}$ and $r_3 \leq 1/2$. Therefore $\gamma_2 = \frac{1}{5}\min\{\frac{q}{1+\frac{q^2}{nL}},\frac{n\mu}{q}\}$ satisfies the above inequality and completes the proof.
\end{proof}

We now use the above result to prove the convergence of QHM iterations~(see Sec.~\ref{sec:discretization}).

\begin{tcolorbox}
\QHM*
\end{tcolorbox}
\begin{proof}
First, we show how one can alternatively write QHM as one-line scheme. The original QHM algorithm is reported here for convenience of the reader
\begin{align*}\tag{QHM}
 	\begin{cases}
 		x_{k+1} = x_k - s((1-a)\nabla f(x_k) + a g_{k+1}) \\
 		g_{k+1} = b g_k +  \nabla f(x_k).
 	\end{cases}
 \end{align*}

We replace the second line of QHM into the first one :
  \begin{align}
      x_{k+1} = x_{k} - s(1-a) \nabla f(x_k) - s \cdot b \cdot a \cdot g_k - a s \nabla  f(x_{k}).
  \end{align}
  Using the first iterate we get: 
  \begin{align}
      -(x_{k}-x_{k-1})/s -(1 -\alpha) \nabla f(x_{k-1}) = a g_{k}.
  \end{align}
  Replacing this into the result of first equation, we get: 
  \begin{align}
      x_{k+1} = x_{k} - s(1-a) \nabla f(x_k) + b ( (x_{k}-x_{k-1}) + s(1 -a) \nabla f(x_{k-1})) - a s\nabla f(x_{k}).
  \end{align}
  By rearrangment, we finally obtain 
 \begin{align}
     x_{k+1} = x_{k} + b(x_k - x_{k-1}) - s \nabla f(x_{k}) + s b(1-a) \nabla f(x_{k-1}).
 \end{align}
The above iterates can be viewed as SIE discretization of \ref{eq:our_model} with the following specific choice of parameters (see the single sequence of iterates of SIE in the last section): 
\begin{align}
    m = (1-a)\sqrt{s},\quad n=a, \quad q=\frac{1-b}{\sqrt{s}}. 
\end{align}
Invoking Thm.~\ref{thm:convergence}, we get the convergence rate for \ref{eq: qhm}. More precisely, choosing $b=1-2\sqrt{\mu s}$ we obtain
\begin{align}
    q = 2\sqrt{\mu}.
\end{align}
The above choice of parameters obeys the constraints in Thm.~\ref{thm:convergence}:
\begin{align}
    m\sqrt{s}  = (1-a)s < s \leq \frac{1}{4L}, \quad n  = a \leq (1-a) = \frac{m}{\sqrt{s}}, 
\end{align}
and
\begin{align}
q\sqrt{s} = 2\sqrt{\mu s} \leq 2\sqrt{sL / 9} \leq 1/3,
\end{align}
since we assumed $s \leq 1/(4L)$, $a\leq1/2$ and $L/\mu\geq 9$. The rate --- thanks to Thm.~\ref{thm:convergence} --- is determined by $\gamma_2 = \frac{1}{5}\min\{\frac{n\mu}{q},\frac{q}{1+q^2/(nL)}) \}$. We conclude the proof by showing that $\gamma_2= a\sqrt{\mu}/8$ in the case of QHM. First, one can readily check that $(n\mu)/(5q) = a\mu/(10\sqrt{\mu})$ holds due to the choice of parameters. Second, with some patience, one can check that the following chain of inequality holds: \[\frac{1}{5}\cdot\frac{q}{1+\frac{q^2}{nL}} = \frac{2\sqrt{\mu}}{5+\frac{20\mu}{aL}}\geq \frac{2a\sqrt{\mu}}{5a+20/9} \geq \frac{a\sqrt{\mu}}{10}.\]
\end{proof}

\section{Proofs for Section~\ref{sec:truncation}} \label{sec:disc_error_proofs}
As stated in the main paper, we consider the following discretization errors:
\begin{align*}
     & \Delta_k^{\text{(\ee)}} := \| X(k\sqrt{s}) - x_{k} \|, \quad x_{k} \; \text{obtained by EE} \\ 
     & \Delta_k^{\text{(\sie)}} := \| X(k\sqrt{s}) - x_{k+1} \|, \quad x_{k} \; \text{obtained by SIE}.
\end{align*}
We define $w_k:= x_k$ for EE and $w_k:= x_{k+1}$ for SIE. We compare the error $\Delta_k = \| X(k\sqrt{s}) - w_{k} \|$ for $k=1$ in the next lemma, assuming $\Delta_0=0$ and $v_0 = V(0)$. This is also called local~(or one-step) integration error.
\begin{tcolorbox}
\trunError*
\end{tcolorbox}
\begin{proof}
 	We introduce the notation $X_k:=X(k\sqrt{s})$, $V_k:=V(k\sqrt{s})$. Our problem setting requires $w_0=X_0$ and $v_0=V_0$. For SIE, $w_k=x_{k+1}$ and we begin from Taylor expansion of $X$ as
 	\begin{align}
 		X_{1} - X_0 = \sqrt{s}\dot{X}_0 + s \ddot{X}_0 + O(s^{\sfrac{3}{2}}),
 	\end{align}
 	and therefore
 	\begin{align}
 	X_1 - w_1 = & X_{1} - X_0 - (w_{1} - w_{0}) + X_0 - w_{0} \nonumber\\
 	= & X_{1} - X_0 - (x_{2} - x_{1}) + X_0 - x_1 \nonumber\\
 	= & \sqrt{s}\dot{X}_0 + s\ddot{X}_{0} + m\sqrt{s}\nabla f(x_1) + n\sqrt{s}v_1  + O(s^{\sfrac{3}{2}}) \nonumber\\
 	= & \sqrt{s}\Big(-m\nabla f(X_0) - n V_0\Big) + s\frac{d}{dt}\Big(-m\nabla f(X_0) - n V_0\Big) \nonumber\\
 	& \quad +m\sqrt{s}\nabla f(x_1) + n\sqrt{s}\Big(\sqrt{s}\nabla f(x_1) + (1-q\sqrt{s})v_0 \Big) + O(s^{\sfrac{3}{2}}).
 	\end{align}
 	where in the third equality we used the fact that, by hypothesis, $X_0-x_1 =0$. 
 	And in particular, since $\frac{d\nabla f(X)}{dt}=\nabla^2 f(X)\dot{X}$,
 	\begin{align}
 	s\frac{d}{dt}\Big(-m\nabla f(X_0) - n V_0\Big) = & -sm\nabla^2 f(X_0)\dot{X_0}-sn\dot{V}_0 \nonumber\\
 	= & -sn\nabla f(X_0) + snqV_0 + sm^2\nabla^2f(X_0)\nabla f(X_0) + smn\nabla^2f(X_0)V_0 . 
 	\end{align}
 	Then it holds that
 	\begin{align}
 	X_1 - w_1 = & -(m\sqrt{s}+ns)\Big(\nabla f(X_0)-\nabla f(x_1)\Big) - n\sqrt{s}(V_0-v_0)  + snq(V_0-v_0) + O(s^{\sfrac{3}{2}})  \leq O(s^{\sfrac{3}{2}}) \label{eq:muk}
 	\end{align}
     and $\Delta_1^{\text{(\sie)}}\leq O(s^{\sfrac{3}{2}})$.
    
     We proceed with the EE iterations (remember: $w_k=x_k$). We expand $\Delta_1^{\text{(\ee)}}$ as
 	\begin{align}
 	X_1 - w_1 = & X_{1} - X_0 - (w_{1} - w_0) + X_0 - x_0 + O(s^{\sfrac{3}{2}}) \nonumber\\
 	= & X_{1} - X_0 - (x_{1} - x_0) + X_0 - x_0 + O(s^{\sfrac{3}{2}}) \nonumber\\
 	= & \sqrt{s}\dot{X}_0 + s\ddot{X}_0 + m\sqrt{s}\nabla f(x_0) + n\sqrt{s} v_0 + O(s^{\sfrac{3}{2}}) \nonumber\\
 	= &\sqrt{s}\Big(-m\nabla f(X_0) -nV_k\Big) + m\sqrt{s}\nabla f(x_0) + n\sqrt{s} v_0  \nonumber\\
 	& \quad + sm\Big(m\nabla^2f(X_0) \nabla f(X_0) + n\nabla^2f(X_k)V_0 \Big) \nonumber\\
 	& \quad - sn\Big(\nabla f(X_0) - qV_0\Big) + O(s^{\sfrac{3}{2}}) \nonumber\\
 	= & -m\sqrt{s}\Big(\nabla f(X_0) - \nabla f(x_0)\Big)  -n\sqrt{s}\Big(V_0 - v_0\Big) +  O(s).
 	\end{align}
     Therefore, we conclude that
     $\Delta_1^{\text{(\ee)}} \leq O(s)$.
 \end{proof}
\begin{tcolorbox}
\globalError*
\end{tcolorbox}

\begin{proof}
The proof is based on the following consequence of strong convexity 
\begin{align}
    \mu \| x - x^* \|^2/2 \leq f(x) - f(x^*).
\end{align}
Using the above inequality together with a straightforward application of triangular inequality we complete the proof: 
\begin{align}
    \| X(k \sqrt{s}) - x_k \| & = \| X(k \sqrt{s}) - x^* +x^*  - x_k \| \nonumber\\
    & \leq \| X(k\sqrt{s}) - x^* \| + \| x_k - x^* \| \nonumber\\ 
    & \leq \sqrt{2} \mu^{-\sfrac{1}{2}} \left( \left( f(X(k\sqrt{s})) - f(x^*)\right)^{\sfrac{1}{2}} + \left( f(x_k) - f(x^*)\right)^{\sfrac{1}{2}}\right).
\end{align}
Replacing the convergence results in Thm.~\ref{thm:continuous}, \ref{thm:convergence}, and \ref{thm:convergence_ee} into the the above bound concludes the proof. 
\end{proof}

\end{document}